\documentclass [12pt,leqno] {article}

\usepackage{url,amsfonts}

\def\dist{\mathop{\rm dist}\nolimits}
\def\supp{\mathop{\rm supp}\nolimits}
\def\loc{\mathop{loc}\nolimits}

\newcommand{\qed}{~\hfill~$\fbox{}$ }
\newcommand{\proces}{( X_t,P^x)} 
\newcommand{\RR}{\mathbb{R}}
\newcommand{\Z}{\mathbb{Z}}
\newcommand{\R}{ \mathbb{R}^{d}}
\newcommand{\N}{\mathbb{N}}
\newcommand{\V}{\mathbb{V}}
\newcommand{\czc}[1]{[ {#1} ]}
\newcommand{\czu}[1]{\{ {#1} \}}

\newcommand{\czcd}[1]{[ {#1} ]^-  }
\newcommand{\dowod}[1][]{{\em Proof#1}.\/ }
\newcommand{\indyk}[1]{{\bf 1}_{#1}}
\newcommand{\scalp}[2]{#1\cdot#2}
\newcommand{\sfera}{ \mathbb{S}}

\newcommand{\Fourier}{ {\cal F}}
\newcommand{\Dynkin}[2]{ {\cal U} #1 (#2) }
\newcommand{\Dynkinr}[2]{ {\cal U}_r #1 (#2) }

\newcommand{\gener}{{\cal A}}
 \def\dist{\mathop{\rm
    dist}\nolimits} 

\newtheorem{lemat}{\indent\sc Lemma}
\newtheorem{prop}{\indent\sc Proposition}
\newtheorem{twierdzenie}{\indent\sc Theorem}
\newtheorem{wniosek}[lemat]{\indent\sc Corollary}
\newtheorem{definicja}{\indent\sc Definition}

\newcounter{conum} 

\begin{document}

\title{Regularity of harmonic functions for anisotropic fractional Laplacian}
\author{Pawe{\l} Sztonyk}
 \footnotetext{ Institute of Mathematics and Computer Science,
  Wroc{\l}aw University of Technology,
   Wybrze{\.z}e Wyspia{\'n}\-skie\-go 27,
   50-370 Wroc{\l}aw, Poland\\
   {\rm e-mail: sztonyk@pwr.wroc.pl
}}
\date{June 4, 2007}
\maketitle

\begin{center}
  Abstract
\end{center}
\begin{scriptsize}
We prove that bounded harmonic functions of anisotropic fractional Laplacians are H\"older continuous under mild regularity assumptions on the corresponding L\'evy measure. Under some stronger assumptions the Green function, Poisson kernel and the harmonic functions are even differentiable of order up to three.
\end{scriptsize}

\footnotetext{2000 {\it MS Classification}:
Primary 47D03, 31C05; Secondary 60J35, 60G51.
Research supported by KBN 1 P03A 026 29 and The Alexander von Humboldt Foundation}
\section{Introduction}\label{s:m}
\setcounter{equation}{0}

Let $d\geq 2$ be a natural number and let $\alpha\in (0,2)$. We investigate a convolution semigroup
of nondegenerate symmetric $\alpha$-stable probability measures $\{P_t\,,\;t> 0\}$ on $\R$ with the corresponding
L\'evy measure $\nu$ (for definitions see Section~\ref{s:p}).
Each $P_t$ has a smooth density $p_t$. 
We consider the {\it potential measure} 
$\V=\int_0^\infty P_t \,dt$
and the {\it potential kernel} 
$$  
  V(x)=\int_0^\infty p_t(x)dt\,,\quad x\in \R\,.
$$
We recall that $V(x)=|x|^{\alpha-d}V(x/|x|)$, but the function may be infinite in some
directions (\cite[pp. 148-149]{BS0}). 
It is known that if $\nu$ is a $\gamma$-measure 
on $\sfera$ (see (\ref{nu_gamma})) with $\gamma>d-2\alpha$, then $V$ is continuous 
hence bounded on $\sfera$ (see \cite[Theorem 1]{BS}).
Conversely, if $V$ is bounded on $\sfera$ then $\nu$ is a $(d-2\alpha)$ - measure on $\sfera$ (\cite[Theorem 2]{BS}).

In the present paper we study H\"older continuity and differentiability of $V$. In what
follows we assume that $\nu$ is a $\gamma$-measure on $\sfera$, see (\ref{nu_gamma}). We denote
$$
  \kappa_0=\gamma-(d-2\alpha).
$$
We note that $\gamma\leq d$, therefore $\kappa_0\leq 2\alpha<4$. The following theorem is our main result on the regularity of the potential kernel (for the definition of the relevant function spaces see below).

\begin{twierdzenie}\label{jadropotencjalu}
  Assume that $\kappa_0>0$. If $\kappa_0\not\in\N$ then the potential kernel $V(x)$ 
  belongs to the local H\"older space $C_{\loc}^{\kappa_0}(\R\setminus\{0\})$. If $\kappa_0\in\N$ then 
  $V(x)\in C_{\loc}^{\kappa_0-}(\R\setminus\{0\})$.
  For every multiindex $\beta$ such that $0\leq |\beta|< \kappa_0$ there exists a constant
  $c=c(\beta)$ such that
  \begin{equation}\label{pochodna_osz}
    |D^{\beta} V(x)|\leq c |x|^{\alpha-d-|\beta|}\,, \quad |x|>0.
  \end{equation}
\end{twierdzenie}

We define operator $\gener$ on smooth 
functions $\varphi$ with compact support in $\R$, $\varphi \in C^\infty_c(\R)$, by letting
\begin{eqnarray}
  \gener\varphi(x) &=&
  \int\limits_{\R}\left(\varphi(x+y)-\varphi(x)-
\scalp{y}{\nabla \varphi(x)}\;\indyk{|y|<1}
\right)\,
  \nu(dy)\nonumber \\
&=& \lim_{\varepsilon\rightarrow 0^+} 
\int\limits_{|y|>\varepsilon}
\left(\varphi(x+y)-\varphi(x)\right)\,
  \nu(dy)\nonumber\,.
\end{eqnarray}
$\gener$ is called 
anisotropic fractional Laplacian and it is a restriction of the infinitesimal generator of
$\{P_t\}$ \cite[Example 4.1.12]{Jc1}. 
In the special case of
$\nu(dy)=c|y|^{-d-\alpha}dy$ we obtain the fractional Laplacian
$\Delta^{\alpha/2}$. For properties of $\Delta^{\alpha/2}$ and a discussion of
equivalent definitions of its harmonic functions 
we refer the reader to \cite{BBsm1999,BKK}.

Harmonic functions corresponding to $\gener$, or $\nu$, are defined by
the mean value property with respect to an appropriate family of harmonic measures,
see Section~\ref{s:p}. We denote
$$
  \kappa_1=\gamma-(d-\alpha),
$$
and we observe that $\kappa_1=\kappa_0-\alpha$ and $\kappa_1\leq \alpha <2$.

\begin{twierdzenie}\label{t:harmoniczne}
  If $\kappa_1>0$ then there exists a constant $\rho\in (0,\kappa_1]$ such that every function $u$ bounded on $\R$ and harmonic in $B(0,1)$ belongs to the local H\"older space $C^{\rho}_{\loc}(B(0,1))$ and
  \begin{equation}\label{e:hh}
    \| u \|_{C^{\rho}(B(0,1-r))} \leq c r^{-\rho} \| u\|_\infty, \quad r\in (0,1/2].
  \end{equation}
If $\kappa_1 \leq \alpha<1$ then we can take $\rho=\kappa_1$ and if $\kappa_1>1$ then $\rho\in (1,\kappa_1)$.
\end{twierdzenie}

In particular every bounded function harmonic in $B(0,1)$ is H\"older continuous if $\kappa_1\in(0,1]$ and even differentiable if $\kappa_1>1$. We note that for $\alpha\geq 1$ an explicit value of $\rho$ can be obtained from Lemma \ref{harmon_holder} and
(\ref{harm_pholder}).

Here are a few comments about the general context of our study. A primary source 
of motivation for this study is the potential theory of the second order elliptic partial differential operators. H\"older continuity of harmonic functions of the second order elliptic partial differential operators
in divergence form was proved independently in \cite{DG} and \cite{Nash}. J. Moser in \cite{Mos} gave another method of the proof using Harnack's inequality and an iteration technique that has been later used successfully in many other situations. In these 
papers the local character of the operators is crucial. Perturbations of strongly elliptic operators
by 'smaller' nonlocal operators can also be dealt with as shown in \cite{MP,HY,Kas}.

Harnack's inequality and H\"older continuity for bounded harmonic functions of purely nonlocal integro-differential operators with kernels $\nu(x,dy)$ were obtained in \cite{BL}. 
We note that the results of \cite{BL} were restricted to absolutely continuous $\nu(x,dy)$ with densities comparable to $|y-x|^{-d-\alpha}$. The class of considered operators then gradually extended, see
\cite{SV} and \cite{SU}. In \cite{BK} the H\"older continuity is proved without the assumption that the density exists.
Instead, the measures $\nu(x,\cdot)$ considered in \cite{BK} satisfy some regularity conditions and estimates from below.
Similar, slightly weaker assumptions yield in \cite{HK} logarithmic modulus of continuity  for harmonic functions: $|u(y)-u(x)|\leq c\|u\|_\infty|\log|y-x||^{-\rho}$.

In our study, initiated in \cite{BSS1} (see also \cite{BSS2,BS0,BS}), we restrict ourselves to the translation invariant case $\nu(x,A)=\nu(0,A-x)=\nu(A-x)$ but we do not generally rely on
the existence of the density of the L\'evy measure $\nu$, and we do not assume 
estimates for $\nu$ from below. We thus allow for much more anisotropy than \cite{BL,SV,SU,BK,HK}. It should be made clear that while our study is analogous to the study of 
symmetric second order elliptic differential operators with mere constant coefficients the variety of different behaviours of the resulting fractional Laplacians makes the
subject very complex and interesting. We believe that the study paves a way to further
generalizations in spirit of \cite{BL} and may also contribute to the general potential theory.

The main focus of the present paper are estimates of the potential kernel $V(x)$ 
(Theorem \ref{jadropotencjalu}). This is a prominent harmonic function of $\gener$ on
$\R\setminus\{0\}$, see \cite{BS0}, and results for more general harmonic functions
(Theorem \ref{t:harmoniczne}) may be considered as consequences of such estimates.
For sufficiently regular
$\nu$ and $\alpha>1/2$ we obtain the differentiability of harmonic functions which is a new result even for
absolutely continuous measures. Our methods are mainly analytic although we use a
probabilistic framework to define some of the objects in question. We do not use Moser's method here because the Harnack's inequality may not hold even if $\nu$ has a bounded density (see \cite{BS0}). We base our development on certain estimates of the derivatives of transition densities (Lemma \ref{s:tdapk}), extending some of the results of
\cite{Lewand,BS,W,Pi,Hi}.

Apart from Theorem \ref{jadropotencjalu} and \ref{t:harmoniczne} we obtain here several results
concerning regularity of transition densities, Green function and Poisson kernel 
for the ball (see Lemma \ref{s:tdapk}, Proposition \ref{Green_regular}, Lemma \ref{Poissonszacbrzeg}, \ref {Poisson_regularnosc} and \ref{l:DP}). These objects were investigated 
in detail in the isotropic case $\nu(dy)=c|y|^{-d-\alpha}dy$. The results obtained include 
estimates of the Green function and Poisson kernel (see \cite{K1,Jk,Mich1,Mich2,BD}), 
estimates of the derivatives of transition density and resolvent (see \cite{BJ,Lewand}) 
and gradient estimates of harmonic functions (see \cite{BKN}). The main difficulty 
of the present study in comparison to the isotropic case is the lack of explicit formulas,
e.g. for the Poisson kernel or the Green function for a ball.
 
The paper is organized as follows. Definitions and preliminaries are
given in Section~\ref{s:p}.
In Section~\ref{s:tdapk} we estimate the derivatives of the transition densities extending
the estimates of the semigroup and potential kernel given in \cite{BS}. At the end of Section~\ref{s:tdapk} 
we prove Theorem \ref{jadropotencjalu}. The continuity estimates of the Green function and its 
derivative are proved in Section~\ref{s:Gf}. In Section~\ref{s:Pk} we investigate the regularity
of the Poisson kernel. Theorem~\ref{t:harmoniczne} is proved in Section~\ref{s:hf}.

\section{Preliminaries}\label{s:p}
The general references for this section are \cite{Ch, Ch1,ChZ,Sato,Br,BGB}.
Recall that $d\geq 2$ and $\alpha\in (0,2)$. We explicitly exclude $d=1$ as
it leads to the well-know rotation invariant case (for which see, e.g., \cite{BBsm1999,BKK}).

For $x\in\R$ and $r>0$ we let $|x|=\sqrt{\sum_{i=1}^d x_i^2}$ and 
$B(x,r)=\{ y\in\R :\: |y-x|<r \}$. We denote $\sfera=\{ x\in\R:\: |x|=1 \}$.
All the sets, functions and measures considered in the sequel will be Borel.
For a measure $\lambda$ on $\R$, $|\lambda|$ denotes
its total mass.
We call $\lambda$ degenerate if there is a proper
linear subspace $M$ of $\R$ such that $\supp(\lambda)\subset M$; 
otherwise we call $\lambda$ {\it nondegenerate}.

In what follows we will consider measures $\mu$ concentrated on $\sfera$.
We will assume that $\mu$ is positive,  finite, nondegenerate (in
particular $\mu\neq 0$), and symmetric:
$$
\mu(D)=\mu(-D)\,,\quad D\subset \sfera\,.
$$
We will call $\mu$ the {\it spectral measure}.
We let
\begin{equation}\label{e:mn}
  \nu(D) = \int_{\sfera} \int_0^\infty \indyk{D}(r\xi) r^{-1-\alpha} \,dr \mu(d\xi)\,,\quad
  D\subset \R\,,
\end{equation}
where $\indyk{D}$ is the indicator function of $D$. 
Note that $\nu$ is symmetric. It is a L\'evy measure on $\R$, i.e.
$$
\int_{\R}\min(|y|^2,\,1)\,\nu(dy)<\infty\,.
$$

For $r>0$ and a function $\varphi$ on $\R$ we consider its dilation
$\varphi_r(y) = \varphi(y/r)$, and we note that $\nu(\varphi_r)=r^{-\alpha} \nu(\varphi)$.
In particular $\nu$ is homogeneous: $\nu(rB)=r^{-\alpha}\nu(B)$ for $B\subset \R$.

We consider, after \cite{BS}, an auxiliary scale of smoothness for $\nu$.
\begin{definicja}
We say that $\nu$ is a $\gamma$-measure on $\sfera$ if there exists a constant $c$ such that
\begin{equation}\label{nu_gamma} 
  \nu(B(x,r)) \leq c r^{\gamma}\,,\quad |x|=1\,,\; 0<r<1/2\,.
\end{equation}
\end{definicja}
Since $\nu(drd\theta)=r^{-1-\alpha}dr\mu(d\theta)$, it is at least a $1$-measure
and at most a $d$-measure on $\sfera$.
If $\nu$ is a $\gamma$-measure with $\gamma>1$, then $\mu$ has no atoms.
Finally, $\nu$ is a $d$-measure if and only if it is absolutely continuous with
respect to the Lebesgue measure and its density function is locally bounded
on $\R\setminus\{0\}$.

We consider the semigroup of {\it stable} probability measures $\{ P_t, t> 0 \}$ with the Fourier transform $\Fourier(P_t)(u)=\exp(-t\Phi(u))$, where
\begin{eqnarray*}
  \Phi(u)
  &  =  & -\int \left(e^{i\scalp{u}{y}}-1-i\scalp{u}{y}\indyk{B(0,1)}(y)\right)\nu(dy) \\
   &  =  & -\int \left(\cos(\scalp{u}{y})-1\right)\nu(dy)
=\frac{\pi}{2\sin\frac{\pi\alpha}{2}\Gamma(1+\alpha)}
  \int_{\sfera}|\scalp{u}{\xi}|^\alpha\mu(d\xi)\,.
\end{eqnarray*}
Here $\scalp{u}{\xi}$ denotes the usual inner product of $u,\xi\in\R$.
Since $\mu$ is finite, homogeneous and nondegenerate,
\begin{equation}\label{Phi}
\Phi(u)=|u|^\alpha \Phi(u/|u|)\approx |u|^{\alpha}\,.
\end{equation}
We also let $\bar{\nu}=\indyk{B(0,1)^c}\nu$, and 
$\tilde{\nu}=\indyk{B(0,1)}\nu$ and consider the corresponding semigroups of measures $\{\tilde{P}_t,\; t\geq 0\}$ and $\{\bar{P}_t,\; t\geq 0\}$ such that
$$
  \Fourier(\bar{P}_t)(u)=
  \int e^{iuy} \bar{P}_t (dy) =
  \exp \left(t\int (\cos(\scalp{u}{y})-1)\bar{\nu}(dy)\right),
  \quad u\in\R\, ,
$$
\begin{equation}\label{FPtilde}
  \Fourier(\tilde{P}_t)(u) =
  \exp\left(t \int (\cos(\scalp{u}{y})-1)
  \tilde{\nu}(dy)\right)\, ,
  \quad u\in\R\, .
\end{equation}

We then have that
\begin{displaymath}
P_t=\tilde{P}_t \ast \bar{P}_t\,,\quad t\geq 0.
\end{displaymath}
The measures $P_t$ and $\tilde{P}_t$
have rapidly decreasing Fourier transform hence they are absolutely continuous with
smooth bounded densities denoted $p_t(x)$
and $\tilde{p}_t(x)$, respectively.
Of course,
\begin{equation}\label{e:wpt}
p_t=
\tilde{p}_t
*
\bar{P}_t
\,, \quad t>0.
\end{equation}

By using (\ref{Phi}) we obtain the {\it scaling property} of $\{p_t\}$:
\begin{equation}\label{e:sc}
  p_t(x)=k^dp_{k^\alpha t}(kx)\,,\quad x\in\R,\, k>0,
\end{equation}
or
\begin{displaymath}
p_t(x)=t^{-d/\alpha}p_1(t^{-1/\alpha}x)\,,\quad x\in \R\, .
\end{displaymath}
In particular,
\begin{equation}\label{ogr_p}
  p_t(x)\leq c t^{-d/\alpha}\,.
\end{equation}
If $\nu$ is a $\gamma$ - measure on $\sfera$ then we have
\begin{equation}\label{starszac}
  p_1(x)\leq c (1+|x|)^{-\alpha-\gamma}\,,
\end{equation}
see \cite{BS} or \cite{W}.
We define the {\it potential measure} of the semigroup $\{ P_t\}$:
$$
  \V(D)=\int_0^\infty P_t(D) dt\,,\quad D\subset \R\,.
$$
Since $\alpha<d$, by (\ref{ogr_p}), $\V$ is finite on bounded subsets of $\R$. 
Let
\begin{equation}\label{eq:pkkk}
  V(x)=\int_0^\infty p_t(x)dt\,,\quad x\in\R\,,
\end{equation}
so that
$$
 \V(D)=\int_D V(x)dx\,,\;\quad D\subset \R\,.
$$
We call $V(x)$ the {\it potential kernel} of the stable semigroup.
By (\ref{e:sc}) 
\begin{equation}
  \label{e:sv}
  V(x)=|x|^{\alpha-d} V\left(x/|x|\right)\,,\quad x\neq 0\, .
\end{equation}
If $\nu$ is a $\gamma$ - measure on $\sfera$ and
$\gamma>d-2\alpha$ then $V(x)$ is continuous on $\sfera$ (\cite[Theorem 1]{BS}).

The L\' evy measure $\nu$ determines the symmetric stable L\' evy process $\proces$
with generating triplet $(0,\nu,0)$. The transition probabilities of the process $\proces$ are 
$P(t,x,A)=P_t(A-x)$, $t>0$, $x\in\R$, $A\subset\R$, and 
$P(0,x,A)=\indyk{A}(x)$. The trajectories are right continuous with left limits and the process is strong Markov with respect 
to the so-called {\it standard filtration} (see, e.g., \cite{Sato}, \cite{Br}, or \cite{BGB}).

For open $U\subset \R$ we denote $\tau_U=\inf\{t \ge 0;\: X_t\not\in U\}$, the
{\it first exit time} of $U$.
We write $\omega_D^x$ for the harmonic measure of (open) $D$:
$$
  \omega_D^x(A)=P^x(\tau_D<\infty\,,\;X_{\tau_D}\in A)\,,\quad x\in \R\,,\; A\subset \R\,.
$$
By the strong Markov property
\begin{equation}
  \label{e:www}
  \omega^x_D(A)=\int \omega^y_D(A)\omega^x_U(dy)\,,\quad\mbox{ if }\; U\subset D\,.
\end{equation}

We say that a function $u$ on $\R$ is
{\it harmonic} in  open  $D\subset \R$ if
\begin{equation}\label{wyrazenie:dH}
  u(x)=E^x u(X_{\tau_U})=\int_{U^c} u(y)\,\omega_U^x(dy)\,,\quad x\in \R,
\end{equation}
for every bounded open set $U$ with the closure $\bar{U}$ contained in
$D$.  It is called {\it regular harmonic} in $D$ if
(\ref{wyrazenie:dH}) holds for $U=D$. If $D$ is unbounded then
$E^x u(X_{\tau_D}) = E^x[\tau_D < \infty\,;\; u(X_{\tau_D})]$
by a  convention.  Under (\ref{wyrazenie:dH}) it will be only assumed that
the expectation in (\ref{wyrazenie:dH}) is well defined (but not necessarily finite).
Regular harmonicity implies harmonicity, and is inherited by subsets $U\subset D$. 
This follows from (\ref{e:www}).

We denote by $p_t^D(x,v)$ the transition density of the process killed 
at the first exit from $D$:
$$p_t^D(x,v)=p(t,x,v)-E^x[\tau_D<t\,;\;p(t-\tau_D,X_{\tau_D},v)]\,,\quad
t>0\,,\; x,v \in \R\,.$$ 
Here $p(t,x,v)=p_t(v-x)$. 
We will assume in the sequel that $D$ is {\it regular}:
$P^x[\inf\{t>0\,:\;X_t\notin D\}=0]=1$ for $x\in D^c$ (see \cite{Ch1, Ch}).
It it well-known that $p_t^D$ is symmetric: $p_t^D(x,v)=p_t^D(v,x)$, $x, v\in \R$ (see, e.g., \cite{ChZ}). 
The strong Markov property yields
\begin{equation}
  \label{e:pt}
   p(t,x,v) = E^x[p(t-\tau_D,X_{\tau_D},v)\,;\;\tau_D<t]\,,
\quad x\in D\,,\; v\in D^c\, ,
\end{equation}
thus $p_t^D(x,v)=0$ if $x\in D$, $v\in D^c$.
We let
  \begin{displaymath}
G_D(x,v)=\int\limits_0^\infty p_t^D(x,v) dt\,,
  \end{displaymath}
and we call $G_D(x,v)$ the {\it Green function} for $D$.  
If $V$ is continuous on $\R\setminus \{0\}$, so that $V(x)\leq
c|x|^{\alpha-d}$, then the strong Markov property yields for $x,v\in D$
\begin{equation}\label{wzornaGreena}
  G_D(x,v)=
V(x,v)-E^xV(X_{\tau_D},v)
=V(x,v)-\int_{D^c}V(z,v)\,\omega_D^x(dz)\,.
\end{equation}
Here $V(x,v)=V(v-x)$. This holds, e.g., if $\nu$ is a $\gamma$-measure on $\sfera$ and $\gamma>d-2\alpha$ (see \cite[Theorem 1]{BS}).
In particular
\begin{equation}\label{Green_podstawa}
 G_D(x,v) \leq c |v-x|^{\alpha-d}\,, \quad x,v\in D.
\end{equation}
The Green function is symmetric: $G_D(x,v)=G_D(v,x)$, continuous
in $D\times D\setminus\{(x,v): x=v\}$, and it vanishes if $x\in D^c$ or $v\in D^c$.

For each $v\in\R$, $V(x,v)$ is harmonic in $x\in\R\setminus\{ v \}$. 
Indeed, if $x\in D$ and $\dist(D,v)>0$ then by (\ref{e:pt}) 
\begin{eqnarray*}
  V(x,v)
  &  =   & \int\limits_0^\infty E^x[p(t-\tau_D,X_{\tau_D},v);\tau_D<t] dt 
  =   E^x V(X_{\tau_D},v)\,.
\end{eqnarray*}
Similarly, the Green function $x\mapsto G_D(x,v)$ is harmonic in $D\setminus\{ v \}$.

We like to note that $V(x,v)$ may be infinite at some points $x\neq v$, see Introduction. 
Thus harmonic functions need to be defined by the mean value property (\ref{wyrazenie:dH})
rather then in terms of the anisotropic fractional Laplacian $\gener$. However,
if a function $\varphi$ belongs to the domain of the infinitesimal generator $\gener'$ of the process,
which is an extension of $\gener$, and $\gener'\varphi(x)=0$ for every $x\in D$ then $\varphi$ is harmonic in $D$ in the sense of (\ref{wyrazenie:dH}). It follows from the fact that $$\varphi(X_t)-\varphi(x)-\int_0^t \gener'\varphi(X_s)ds$$ is
a martingale with respect to $P^x$ (see, e.g., \cite{EtK}).

By Ikeda--Watanabe formula \cite{IW} we have
\begin{equation}\label{IW}
  \omega_D^x(A) = \int_D G_D(x,v)\nu(A-v)dv\,,\quad \mbox{ if }\; \dist(A,D)>0\,.
\end{equation}
If the boundary of $D$ is smooth or even Lipschitz 
then
$$
  \omega_{D}^x(\partial D)=0\,,\quad x\in D\,,
$$
see (\cite{Sztonyk}, \cite{Millar}, \cite{Wu}). 
In this case $\omega_D^x$ is absolutely continuous with respect to
the Lebesgue measure on $D^c$.
Its density function, or the Poisson kernel, is given by the formula (see \cite{BS})
\begin{equation}
  \label{e:P}
  P_D(x,z)=\int_{z-D} G_D(x,z-v) \nu(dv)\,,\quad x\in D\,.
\end{equation}
Note that such $D$ is regular, which follows from scaling and the fact that
\begin{equation}
  p_t(x)>0\,,\quad x\in\R\quad (t>0),
\end{equation}
(see \cite{Taylor} or \cite[Lemma 5]{Pi}).
In particular the above considerations apply to $D=B(0,1)$.
 
It follows from (\ref{Phi}) that for every $r>0$ and $x\in \R$ the
$P^x$ distribution of $\{X_t\,,\,t\geq 0\}$ is the same as the $P^{rx}$
distribution of $\{r^{-1}X_{r^{\alpha}t}\,,\,t\geq 0\}$. 
In consequence,
\begin{equation}
  \label{e:scx}
\omega^x_D(A)=\omega^{rx}_{rD}(rA)\, ,
\end{equation}
which we will call scaling, too.
It yields that for $u$ harmonic on $D$, the dilation,
$u_r$, is harmonic on $rD$. A similar remark applies to translations. 

Let $\N_0=\{0,1,2,\dots\}$. For a multiindex $\beta=(\beta_1,\dots,\beta_d)\in\N_0^d$ we denote
$|\beta|=\beta_1+\dots+\beta_d$. 
Let $s\geq 0$ and $s=\czc{s}+\czu{s}$ where $\czc{s}\in\Z$ and
$\czu{s}\in[0,1)$. For every $s>0$, $s\not\in\N$ and open set $D\subset\R$ we denote by $C^s(D)$ 
the H\"older space of order $s$ (see, e.g., \cite{Trib}), i.e.,
$$
  C^s(D)  =  \{f\in C^{\czc{s}}(D): \|f\|_{C^s(D)}<\infty\} 
$$
where $C^0(D)$ denotes the set of bounded and uniformly continuous functions on $D$ and for $n\in\N$, $ C^n(D)=\{f:\: D^\beta f\in C^0(D) \mbox{ for all } |\beta|\leq n  \}$, $
  \|f\|_{C^n(D)}=\sum_{|\beta|\leq n}\|\indyk{D}\cdot D^\beta f\|_\infty
$, and 
$$
\|f\|_{C^s(D)}=\|f\|_{C^{\czc{s}}(D)}+
           \sum_{|\beta|=\czc{s}} \sup_{x,y\in D}\frac{|D^\beta f(y)-D^\beta f(x)|}{|y-x|^{\czu{s}}}.
$$
We write $f\in C^{s-}(D)$ if $f\in C^{s-\delta}$ for every $\delta\in (0,s]$.
By $C^s_{\loc}(D)$ we denote the local H\"older space, i.e. a function $f$ belongs
to $C^s_{\loc}(D)$ if $f\in C^s(U)$ for every bounded open set $U$ with the closure $\bar{U}$ contained in $D$.

We use $c$ (with subscripts) to denote finite positive constants
which depend only on the measure $\mu$ and the constant $\gamma$, the dimension $d$ and
the index
$\alpha$. Any {\it additional} dependence
is explicitly indicated by writing, e.g., $c=c(\beta,n)$.
The value of $c$, when used without subscripts, may
change from place to place.

\section{Transition density and potential kernel}\label{s:tdapk}
We first investigate the regularity of the density of the semigroup
$\{\tilde{P}_t,t\geq 0\}$ generated by the measure $\tilde{\nu}=\indyk{B(0,1)}\nu$.
Let $S(\R)$ denote the usual Schwarz space of smooth rapidly decreasing functions 
(see, e.g., \cite{Rudin}).
\begin{lemat}\label{Schwarz}
$\Fourier(\tilde{P}_t)\in S(\R)$.
\end{lemat}
\dowod Let $\varphi(u,y)=\cos(\scalp{u}{y})-1$, and $\tilde{\Phi}(u)= \int \varphi(u,y)\tilde{\nu}(dy)$, see (\ref{FPtilde}),
  where $u=(u_1,...,u_d)$, $y=(y_1,...,y_d)$. For every $\beta=(\beta_1,...,\beta_d)\in\N_0^d$, $|\beta|>0$, we have
  $$
    D^\beta_u \varphi (u,y) = y^\beta \cos^{(|\beta|)}(\scalp{u}{y})
  $$
  where $y^\beta=y_1^{\beta_1}\cdot ... \cdot y_d^{\beta_d}$
  and $\cos^{(n)}(\theta)$ denotes the n-th derivative of the function $\RR\ni \theta\to \cos(\theta)$. Hence
  $$
    |D^\beta_u \varphi (u,y)| \leq |y|^{|\beta|+1}(|u|\wedge\frac{1}{|y|})\,,\quad |\beta|=1,3,5,\dots,
  $$
  and
  $$
    |D^\beta_u \varphi (u,y)| \leq |y|^{|\beta|}\,,\quad |\beta|=2,4,6,\dots .
  $$
  Therefore, we can change the order of differentiation and integration, obtaining
  $$
    D^\beta \tilde{\Phi}(u) = \int D^\beta_u \varphi (u,y)\tilde{\nu}(dy).
  $$
  Hence
  $$
    |D^\beta \tilde{\Phi}(u)| \leq c(\beta)g(u),\quad 
    u\in\R,\, |\beta|=1,2,3,4,\dots,
  $$ 
where
  $$
g(u)=
\left\{
\begin{array}{lcl}
   1  & \mbox{for }  & |\beta|=2,4,6,\dots,\\
  |u|\wedge 1& \mbox{for }  & |\beta|=1,3,5,\dots \mbox{  and  } |\beta|>\alpha,\\
   |u|\wedge |u|^{\alpha-1} & \mbox{for }  & |\beta|=1 \mbox{  and  } \alpha>1,\\
   |u|\wedge (1+|\log|u||) & \mbox{for } & |\beta|=\alpha=1.
\end{array}\right.
$$
  We use the multivariable version of Faa di Bruno's formula (see \cite{H}) for the function 
  $f(u)=\Fourier(\tilde{P}_t)(u)=\exp (t\tilde{\Phi}(u))$ to obtain for every $k\in\N$
  $$  
 \frac{\partial^k }{\partial u_1\dots \partial u_k}f(u) = f(u) \sum_{\pi\in \cal{P}} 
 \prod_{A\in\pi}\frac { \partial^{|A|} (t\tilde{\Phi}) }{\prod_{j\in A} \partial u_j }(u),
  $$
  where $\cal{P}$ denotes the set of all partitions of the set $\{1,\dots,k \}$.
  The formula is still valid if some of the variables $u_1,\dots,u_k$ 
  on the left hand side denote the same 
  variable and then the corresponding terms in the sum are multiply counted.
  Since $g(u)\leq 1+|u|$, we get
  $$
    |D^\beta f(u)|\leq c(\beta,t) (1+|u|)^{|\beta|}|f(u)| \leq 
    c(\beta,t)(1+|u|)^{|\beta|} e^{-c|u|^\alpha}
  $$
\qed

Lemma \ref{Schwarz} yields that for every $\beta\in\N_0^d$ and $q>0$ there exists
a constant $c=c(\beta,q)$ such that
\begin{equation}\label{tildaszac}
  |D^\beta \tilde{p}_1 (y)| \leq c (1+|y|)^{-q}\,,\quad y\in\R.
\end{equation}

The proof of the following lemma is based on the ideas of the proofs of \cite[Theorem 3]{Pi} and
\cite[Lemma 6]{BS}. In what follows we use the convention $D^0={\rm Id}$.

\begin{lemat}\label{oszacowanie}
If $\nu$ is a $\gamma$-measure then for every $\beta\in\N_0^d$ there exists a constant $c=c(\beta)$ such that
$$
  |D^\beta p_1(y)|\leq c (1+|y|)^{-\alpha-\gamma}\,, \quad y\in\R.
$$
\end{lemat}
\dowod We have
$$
p_1(y)=\int \tilde{p}_1(y-z) \bar{P}_1(dz).
$$
Also,
$$
D^\beta p_1(y) = \int D^\beta \tilde{p}_1(y-z) \bar{P}_1(dz),
$$ 
because $\bar{P}_1$ is finite and $D^\beta \tilde{p}_1$ is bounded and continuous.
By \cite[Corollary 4]{BS} we have
$$
  \bar{P}_1(B(y,\lambda))\leq c \lambda^\gamma(1+\lambda^\alpha)|y|^{-\alpha-\gamma}, \quad y\in \R,\, \lambda>0.
$$
Hence for $q= \gamma+\alpha +1$ by (\ref{tildaszac}) we obtain
\begin{eqnarray*}
  |D^\beta p_1(y)|  &  =   & | \int_{\R} D^\beta \tilde{p}_1(y-z)\bar{P}_1(dz)| \\
                    & \leq & \int_{\R} |D^\beta \tilde{p}_1(y-z)|\bar{P}_1(dz)\\
                    & \leq & \int_{\R} c(\beta) (1+|y-z|)^{-q} \bar{P}_1(dz)\\
                    &  =   & c(\beta) \int_0^1 \bar{P}_1(\{z\,:\; (1+|z-y|)^{-q}>s\})ds \\
                    & \leq & c(\beta) \int_0^1 \bar{P}_1(B(y, s^{-1/q}))ds\\
                    & \leq & c(\beta) \int_0^{1}(s^{-1/q})^\gamma(1+(s^{-1/q})^\alpha)
                             |y|^{-\gamma-\alpha}ds \\
                    &  =   & c(\beta) |y|^{-\gamma-\alpha} 
                             \left[\int_0^1 s^{-\gamma/q}ds +
                                   \int_0^1 s^{-(\gamma+\alpha)/q}ds
                             \right] \\
                    &  =   & c(\beta) |y|^{-\gamma-\alpha}.
\end{eqnarray*}
Since $D^\beta p_1$ is bounded we also have 
$ |D^\beta p_1(y)| \leq c(\beta) (1+|y|)^{-\alpha-\gamma}$.
\qed

We note that explicit estimates for the first derivative of the transition 
density in the isotropic case of $\nu(dy)=c|y|^{-d-\alpha}dy$ are also given in \cite[Lemma 5]{BJ}. In this case we have $|\nabla_y p(1,y)|\leq c|y|(1+|y|)^{-d-2-\alpha}$.

\begin{lemat}\label{holderp}
  For every $\beta\in\N_0^d$ there exists a constant $c=c(\beta)$ such that for every $R>0$ we have
  \begin{equation}
    |D^\beta p_1(y)-D^\beta p_1(x)|  \leq  c (1+R)^{-\alpha-\gamma}(|y-x|\wedge 1),
  \end{equation}
  provided $x,y\in B(0,R)^c\,,\; |y-x|<R/4$.
\end{lemat}
\dowod By the Lagrange's theorem we obtain
$$
  D^\beta p_1(y)-D^\beta p_1(x) = \scalp{\nabla D^\beta p_1(x+\theta(y-x))}{(y-x)},
$$
where $\theta\in [0,1]$. If $|y-x|<1$ we get by Lemma \ref{oszacowanie}
\begin{eqnarray*}
  |D^\beta p_1(y)-D^\beta p_1(x)|  & \leq & |\nabla D^\beta p_1(x+\theta(y-x))|\cdot |y-x| \\
                                   & \leq & c (1+|x+\theta(y-x)|)^{-\alpha-\gamma}|y-x| \\
                                   & \leq & c (1+(3/4)R)^{-\alpha-\gamma}|y-x| \\
                                   & \leq & c (1+R)^{-\alpha-\gamma}|y-x|.
\end{eqnarray*}
For $|y-x|\geq 1$ we have
$$
  |D^\beta p_1(y)-D^\beta p_1(x)| \leq 2 \sup_{z\in B(0,R)^c} |D^\beta p_1(z)| \leq
  c(1+R)^{-\alpha-\gamma}.
$$
\qed

\dowod[ of Theorem \ref{jadropotencjalu}]
We will first prove that for every $\beta\in\N_0^d$ such that $|\beta|<\kappa_0=\gamma-(d-2\alpha)$, we have
\begin{equation}\label{pochodna}
    D^\beta V(x) = \int_0^\infty D^\beta_x p_t(x) dt, \;\; |x|>0.
  \end{equation}
By scaling property (\ref{e:sc}) we have for every $\beta\in\N_0^d$
$$
  D_x^\beta p_t(x)=k^{d+|\beta|}(D_x^\beta p_{k^{\alpha}t})(kx),
  \;\; x\in\R,\; t>0,\; k>0.
$$
For $k=t^{-1/\alpha}$ we get by Lemma \ref{oszacowanie}
\begin{eqnarray*}
  |D^\beta_x p_t(x)| &  =   & t^{\frac{-d-|\beta|}{\alpha}}|(D^\beta p_1)(t^{-1/\alpha}x)| \\
                   & \leq & c(\beta) t^{\frac{-d-|\beta|}{\alpha}}(1+t^{-1/\alpha}|x|)^{-\alpha-\gamma},
  \;\; x\in\R,\; t>0,
\end{eqnarray*}
and this yields (\ref{pochodna}) and the continuity of $D^\beta V(x)$ for $x\neq 0$.
Hence we have
\begin{eqnarray}\label{skal_pv}
  D^\beta V(x)  &  =  & \int_0^\infty D^\beta_x p_t(x) dt \\ \nonumber
                &  =  & k^{d+|\beta|} \int_0^\infty (D_x^\beta p_{k^{\alpha}t})(kx) dt \\
                &  =  & k^{-\alpha+d+|\beta|} \int_0^\infty (D_x^\beta p_{s})(kx) ds 
                        = k^{-\alpha+d+|\beta|}(D^\beta V)(kx),\nonumber
\end{eqnarray}
where $k>0$. This yields (\ref{pochodna_osz}).\\
Let $|y-x|<1/4$ and $|x|>1$, $|y|>1$. By scaling and 
Lemma \ref{holderp}
\begin{eqnarray*}
  |D^\beta p_t(y)-D^\beta p_t(x)|  &    = & t^{\frac{-d-|\beta|}{\alpha}}|(D^\beta p_1)(t^{-1/\alpha}y)-(D^\beta p_1)(t^{-1/\alpha}x)| \\
                   & \leq & t^{\frac{-d-|\beta|}{\alpha}} 
                            c(1+t^{-1/\alpha})^{-\alpha-\gamma}
                            (|t^{-1/\alpha}y-t^{-1/\alpha}x|\wedge 1)\\
                   &   =  & c t^{\frac{-d-|\beta|-1}{\alpha}}
                            (1+t^{-1/\alpha})^{-\alpha-\gamma}
                            (|y-x|\wedge t^{1/\alpha})\,,\quad t>0.
\end{eqnarray*}
This yields
\begin{eqnarray*}
  |D^\beta V(y)-D^\beta V(x)|  & \leq & c \int_0^\infty t^{\frac{-d-|\beta|-1}{\alpha}}
                                        (1+t^{-1/\alpha})^{-\alpha-\gamma}
                                        (|y-x|\wedge t^{1/\alpha})dt \\
                               &   =  & c \int_0^{|y-x|^\alpha} t^{\frac{-d-|\beta|-1}{\alpha}}
                                        (1+t^{-1/\alpha})^{-\alpha-\gamma}t^{1/\alpha} dt \\
               &      & + \; c \int_{|y-x|^\alpha}^\infty t^{\frac{-d-|\beta|-1}{\alpha}}
                            (1+t^{-1/\alpha})^{-\alpha-\gamma}
                            |y-x|dt \\
               &   =  & I \; + \; II.
\end{eqnarray*}
Let $\kappa_0\not\in\N$ and $|\beta|=\czc{\kappa_0}$. We have
$$
  I  \leq  c \int_0^{|y-x|^\alpha} t^{\frac{-d-|\beta|+\alpha+\gamma}{\alpha}}dt 
      =    c |y-x|^{2\alpha+\gamma-d-|\beta|} = c |y-x|^{\czu{\kappa_0}},
$$
\begin{eqnarray*}
  II & \leq & c \int_{|y-x|^\alpha}^1 t^{\frac{-d-|\beta|-1+\alpha+\gamma}{\alpha}}
              |y-x|dt + c \int_1^\infty t^{\frac{-d-|\beta|-1}{\alpha}} |y-x|dt \\
     & \leq & c (|y-x|^{\czu{\kappa_0}}+|y-x|)\leq c |y-x|^{\czu{\kappa_0}}.
\end{eqnarray*}
If $|y-x|\geq \frac{1}{4}$ then also we can write
\begin{eqnarray*}
  |D^\beta V(y)-D^\beta V(x)|  & \leq & 2\cdot 
  \left(\frac{|y-x|}{4}\right)^{\czu{\kappa_0}}\cdot\sup_{z\in B(0,1)^c} |D^\beta V(z)|.
\end{eqnarray*}
Therefore, by (\ref{skal_pv}) we have that for every $r>0$
and $x,y\in B(0,r)^c$,
\begin{eqnarray}\label{szac_holder_pv}
  |D^\beta V(y)-D^\beta V(x)| &   = & r^{\alpha-d-|\beta|}|(D^\beta V)(y/r)-(D^\beta V)(x/r)| \nonumber \\ 
              &\leq & c r^{\alpha-d-|\beta|} (|y-x|/r)^{\czu{\kappa_0}}\, .
\end{eqnarray}
Thus if $\kappa_0\not\in\N$ then $V\in C^{\kappa_0}(B(0,r)^c)$ for every $r>0$, or 
$V\in C^{\kappa_0}_{\loc}(\R\setminus\{ 0 \}).$\\
We now assume that $\kappa_0\in\N$ and $|\beta|=\kappa_0-1$. Let $|x|>1,\; |y|>1,\; |y-x|<1/4$.
We have
\begin{eqnarray*}
  I &   =   &\int_0^{|y-x|^\alpha} t^{\frac{-d-|\beta|-1}{\alpha}}
             (1+t^{-1/\alpha})^{-\alpha-\gamma}t^{1/\alpha} dt \\
    & \leq & c \int_0^{|y-x|^\alpha} t^{\frac{1-\alpha}{\alpha}}dt \\
    &  =   & c |y-x|,
\end{eqnarray*}
and
\begin{eqnarray*}
  II &  =   & \int_{|y-x|^\alpha}^\infty t^{\frac{-d-|\beta|-1}{\alpha}}
                            (1+t^{-1/\alpha})^{-\alpha-\gamma}
                            |y-x|dt \\
     & \leq & c \int_{|y-x|^\alpha}^1 t^{-1}
              |y-x|dt + c \int_1^\infty t^{\frac{-d-|\beta|-1}{\alpha}} |y-x|dt \\
     &  =   & c |y-x|(\log(\frac{1}{|y-x|})+1).
\end{eqnarray*}
This yields
$$
  |D^\beta V(y)-D^\beta V(x)|\leq c |y-x|\log(2+\frac{1}{|y-x|}).
$$
If $|y-x|\geq \frac{1}{4}$ then we can write
\begin{eqnarray*}
  |D^\beta V(y)-D^\beta V(x)|  & \leq & 16 |y-x|\log(2+\frac{1}{|y-x|})\sup_{z\in B(0,1)^c} |D^\beta V(z)|.
\end{eqnarray*}
Thus, for every $r>0$, $x,y\in B(0,r)^c$ by (\ref{skal_pv}) we get 
\begin{eqnarray}\label{logoholder}
  |D^\beta V(y)-D^\beta V(x)| &   = & r^{\alpha-d-|\beta|}|(D^\beta V)(y/r)-(D^\beta V)(x/r)|\nonumber \\   
              &\leq & c r^{\alpha-d-|\beta|} (|y-x|/r) \log(2+\frac{r}{|y-x|}).
\end{eqnarray}
In particular, $V(x)\in C^{\kappa_0-p}(B(0,r)^c)$ for every $p\in (0,1)$.
\qed

We note that if the L\'evy measure $\nu$ is absolutely continuous and has
a locally bounded density, i.e., $\gamma=d$, then $\kappa_0=2\alpha$ and 
the potential kernel is differentiable for $\alpha>1/2$ and $2\alpha$-H\"older
continuous for $\alpha\in (0,1/2)$.

\section{Green function}\label{s:Gf}

In this section we consider the regularity of the Green function of the unit ball.

For every finite measure $m$ we denote by $\V m(x)$ its potential at $x$:
$$
  \V m(x)=\int_{\R} V(x-z) m(dz).
$$ 

In what follows for every $s\in\RR$ we denote
  $$
\czcd{s}=
\left\{
\begin{array}{lcl}
  \czc{s}   & \mbox{for }  & s\not\in\Z,\\
  s-1       & \mbox{for }  & s\in\Z.
\end{array}\right.
$$

\begin{lemat}\label{potencjalmiary} Let $m$ be a finite measure with support $S$ and $\kappa_0=\gamma-(d-2\alpha)>0$.
  If $\kappa_0\not\in\N$ then the potential $\V m(x)$ belongs to the local H\"older space $C_{\loc}^{\kappa_0}(S^c)$. If 
  $\kappa_0\in\N$ then $\V m(x)\in C_{\loc}^{\kappa_0-}(S^c)$. Moreover there is $c$ such that
  for every $r>0$ and $|\beta|=\czcd{\kappa_0}$ we have
  \begin{equation}\label{holder_pm}
    |D^\beta \V m(y)-D^\beta \V m(x)|   \leq   c r^{\alpha-d-|\beta|} f\left(\frac{|y-x|}{r}\right)|m|,
  \end{equation}
  provided $\dist(x,S)>r,\, \dist(y,S)>r$, where
  $$
f(\rho)=
\left\{
\begin{array}{lcl}
  \rho^{\czu{\kappa_0}}   & \mbox{for }  & \kappa_0\not\in\N,\\
  \rho\log(2+\frac{1}{\rho})    & \mbox{for }  & \kappa_0\in\N.
\end{array}\right.
$$
For every $|\beta|<\kappa_0$ we have
  \begin{equation}\label{pochodna_pm_osz}
    |D^\beta \V m(x)| \leq c \dist(x,S)^{\alpha-d-|\beta|}|m|\,,\quad x\in S^c.
  \end{equation}
  
\end{lemat}
Lemma \ref{potencjalmiary} follows easily from Theorem \ref{jadropotencjalu}, (\ref{szac_holder_pv}) and (\ref{logoholder}) and so we omit the proof.

In what follows we denote $B=B(0,1)$.

\begin{prop}\label{Green_regular}
  Let $\kappa_0=\gamma-(d-2\alpha)>0$.
  If $\kappa_0\not\in\N$ then for every $v\in B$ the function $g(x)=G_B(x,v)$ belongs to the local H\"older space $C_{\loc}^{\kappa_0}(B\setminus \{v\})$. If 
  $\kappa_0\in\N$ then $g(x) \in C_{\loc}^{\kappa_0-}(B\setminus\{v\})$. Moreover there is $c$ such that for every $r\in (0,1)$ and $|\beta|=\czcd{\kappa_0}$ we have
  \begin{equation}\label{holder_Green}
    |D^\beta  g(y)-D^\beta  g(x)|   \leq   c r^{\alpha-d-|\beta|} f\left(\frac{|y-x|}{r}\right),
  \end{equation}
  provided $x,y\in B(0,1-r)\cap B(v,r)^c$, where
  $$
f(\rho)=
\left\{
\begin{array}{lcl}
  \rho^{\czu{\kappa_0}}   & \mbox{for }  & \kappa_0\not\in\N,\\
  \rho\log(2+\frac{1}{\rho})    & \mbox{for }  & \kappa_0\in\N.
\end{array}\right.
$$
For every $|\beta|<\kappa_0$ we have
  \begin{equation}\label{pochodna_Green_osz}
    |D^\beta  g(x)| \leq c |v-x|^{\alpha-d-|\beta|}\,,\quad x\neq v,\, |x|<3/4.
  \end{equation}
\end{prop}
\dowod
Since $\supp(\omega^v_B)\subset B^c$ 
the Proposition  follows from (\ref{wzornaGreena}), Theorem \ref{jadropotencjalu}, (\ref{szac_holder_pv}), (\ref{logoholder}) and Lemma \ref{potencjalmiary}.
\qed

\begin{lemat}\label{calka_gamma_miary}
  For every $0<a<\gamma$ there exists a constant $c=c(a)$ such that
  \begin{equation}
    \int_{B(x,\rho)} |z-x|^{-a} \nu(dz) 
      \leq c\rho^{\gamma-a}|x|^{-\alpha-\gamma}\,,\quad x\in\R\setminus\{0\},\, 0<\rho<\frac{|x|}{2}.
  \end{equation}
\end{lemat}
\dowod Let $x\in\R\setminus\{0\}$, $0<\rho<|x|/2$. Then
\begin{eqnarray*}
  \int_{B(x,\rho)} |z-x|^{-a} \nu(dz) &   =  & \sum_{n=0}^\infty 
                                               \int_{B(x,2^{-n}\rho)\setminus B(x,2^{-n-1}\rho)}
                                               |z-x|^{-a}\nu(dz) \\
                                    & \leq & \sum_{n=0}^\infty (2^{-n-1}\rho)^{-a} 
                                             \nu(B(x,2^{-n}\rho)) \\
                                    & \leq & \sum_{n=0}^\infty (2^{-n-1}\rho)^{-a} 
                                             |x|^{-\alpha}\nu(B(x/|x|,2^{-n}\rho/|x|)) \\
                                    & \leq & \sum_{n=0}^\infty (2^{-n-1}\rho)^{-a} 
                                             c |x|^{-\alpha-\gamma}(2^{-n}\rho)^\gamma \\
                                    &    = & \frac{2^a}{1-2^{a-\gamma}}c\rho^{\gamma-a}
                                             |x|^{-\alpha-\gamma}.
\end{eqnarray*}
\qed

We define
$$ 
  s(x)=E^x \tau_{B(0,1)}=\int_{B(0,1)} G(x,v)dv\,.
$$

The following result of M. Lewandowski is consequence of symmetry and nondegeneracy of the spectral measure $\mu$ (for the proof see \cite[Lemma 10]{BS0}).

\begin{lemat}\label{ET_est} There exists 
$c$ such that
$$
  s(x) \leq c (1-|x|^2)^{\alpha/2}\,,\quad |x|<1.
$$
\end{lemat}

\begin{lemat} If $\kappa_1=\gamma-(d-\alpha)>0$ then
\begin{equation}\label{Green_szac}
  G_B(x,v)\leq c(1-|v|)^{\alpha/2}\,,\quad |x|<1/2,\, |v|>3/4.
\end{equation}
\end{lemat}
\dowod We refer to the proof of Lemma 16 in \cite{BS}. The only modification we need to make is Lemma \ref{calka_gamma_miary} above with $a=d-\alpha$ to estimate $U_r g_n(v)$ in $\cite{BS}$.
\qed

\section{Poisson kernel}\label{s:Pk}

\begin{lemat}\label{miara_L}
 There exists a constant $c$ such that for every $r>0$ we have
\begin{equation}
  \nu(B(z,r)) \leq c r^{\gamma}|z|^{-\gamma}(|z|-r)^{-\alpha}\,,\quad |z|>r.
\end{equation}
\end{lemat}
\dowod For $|z|>2r$ we have
$$
  \nu(B(z,r))=|z|^{-\alpha}\nu(B(\frac{z}{|z|},\frac{r}{|z|})) \leq 
  c r^{\gamma}|z|^{-\alpha-\gamma},
$$
whereas for $r<|z|\leq 2r$
$$
  \nu(B(z,r)) \leq \nu (B(0,|z|-r)^c) = c (|z|-r)^{-\alpha}.
$$
\qed

\begin{lemat}\label{Poissonszacbrzeg}
If $\kappa_1=\gamma-(d-\alpha)>0$ then
\begin{equation}\label{Poisson_szac}
  P_B(x,z)\leq c |z|^{-\gamma} (|z|^2-1)^{-\alpha/2}\,,\quad |x|<1/2,\; |z|>1.
\end{equation}
\end{lemat}
\dowod By (\ref{e:P}) and (\ref{Green_podstawa}) we have
$$
  P_B(x,z) = \int_{B(z,1)} G_B(x,z-v)\nu(dv) \leq \int_{B(z-x,3/2)} c|v-(z-x)|^{\alpha-d} \nu(dv).
$$
For $|z|>3$ the estimate (\ref{Poisson_szac}) follows from Lemma \ref{calka_gamma_miary} with $a=d-\alpha$, and Lemma \ref{miara_L}.
For $1<|z|<3$ we have by (\ref{Green_szac}), Lemma \ref{calka_gamma_miary} with $a=d-\alpha$, and Lemma \ref{miara_L},
\begin{eqnarray*}
  P_B(x,z)  &  =   & \int_{B(z,1)} G_B(x,z-v)\nu(dv) \\
            & \leq & \int_{B(z,1)\setminus B(z,3/4)} c(1-|z-v|)^{\alpha/2} \nu(dv)+ 
                     \int_{B(z,3/4)} c|v-(z-x)|^{\alpha-d} \nu(dv) \\
            & \leq & \int_{B(0,|z|-1)^c} |v|^{\alpha/2} \nu(dv) +c \\
            & \leq & c (|z|-1)^{-\alpha/2}.
\end{eqnarray*}
\qed

The main result of this section is the following lemma.

\begin{lemat}\label{Poisson_regularnosc}
If $\kappa_1=\gamma-(d-\alpha)>0$ and $\kappa_1\neq 1$ then for every $|z|>1$ the function $k_z(x)=P_B(x,z)$
belongs to the H\"older space $C^{\kappa_1}(B(0,1/2))$. If $\kappa_1=1$ then
$k_z(x)\in C^{1-}(B(0,1/2))$. Moreover there is a constant $c$ such that
for every $|\beta|=\czcd{\kappa_1}$ we have
\begin{eqnarray}\label{Poisson_holder}
  |D^\beta k_z(y)-D^\beta k_z(x)| & \leq  & cf(|y-x|)
                                |z|^{-\gamma} (|z|-1)^{-\alpha},\\
                      &       & x,y\in B(0,1/2),\; |z|>1,\nonumber
\end{eqnarray}
where
$$
f(\rho)=
\left\{
\begin{array}{lcl}
  \rho^{\czu{\kappa_1}}   & \mbox{for }  & \kappa_1\not\in\N,\\
  \rho\log(\frac{2}{\rho})    & \mbox{for }  & \kappa_1\in\N.
\end{array}\right.
$$
If $\kappa_1>1$ then for every $i\in\{1,\dots,d\}$ we have
\begin{equation}\label{Poisson_pochodna_szac}
  |D_i k_z(x)| \leq 
  c |z|^{-\gamma} (|z|-1)^{-\alpha}\,,\quad |x|<1/2,\; |z|>1.
\end{equation}
\end{lemat}
\dowod We assume first that $\kappa_1=\gamma-(d-\alpha)>1$.
For every fixed $x\in B(0,1/2)$, $i\in\{1,\dots,d \}$ and $|z|>1$ we denote
$$
  f_h(v)=\frac{ G_B(x+he_i,z-v)-G_B(x,z-v) }{h},\,\quad h\in\RR, v\in B,
$$
where $\{e_1,\dots e_d \}$ is the standard orthonormal basis in $\R$.
We will prove that the family of functions 
$\{ f_h(v), |h|<1/16 \}$ is uniformly integrable with respect to the measure
$\indyk{B(z,1)}\nu(dv)$. Thus we need 
to estimate the integral
\begin{eqnarray*}
  \int_{B(z,1)}|f_h(v)|^{1+\varepsilon} \nu(dv) & = & \int_{B(z-x,2|h|)} +
   																							 \int_{B(z-x,1/8)\setminus B(z-x,2|h|)}
   																							+ \int_{B(z,1)\setminus B(z-x,1/8)} \\
   																							& = & A+B+C,
\end{eqnarray*}
where $0<\varepsilon<\frac{\gamma}{d-\alpha+1}-1$.
By Lemma \ref{calka_gamma_miary} with $a=(d-\alpha)(1+\varepsilon)$ 
and (\ref{Green_podstawa}) we have
\begin{eqnarray*}
  A & \leq & \frac{2^{1+\varepsilon}}{|h|^{1+\varepsilon}} 
               [\int_{B(z-x-he_i,3|h|)} (G_B(x+he_i,z-v))^{1+\varepsilon}\nu(dv) \\ 
    &      &      + \int_{B(z-x,3|h|)} (G_B(x,z-v))^{1+\varepsilon}\nu(dv)] \\
    & \leq & \frac{2^{1+\varepsilon}}{|h|^{1+\varepsilon}} [ \int_{B(z-x-he_i,3|h|)} 
                                  c|z-x-he_i-v|^{(\alpha-d)(1+\varepsilon)}\nu(dv) \\
    &      &  + \int_{B(z-x,3|h|)} c |z-x-v|^{(\alpha-d)(1+\varepsilon)}\nu(dv)] \\
    & \leq & \frac{c}{|h|^{1+\varepsilon}} |h|^{\gamma+(\alpha-d)(1+\varepsilon)}
             |z|^{-\alpha-\gamma} \\
    & \leq & c |z|^{-\alpha-\gamma}.
\end{eqnarray*}
Let $L=\czc{-\log_2(8|h|)}$. By (\ref{pochodna_Green_osz}) with $|\beta|=1$ and Lagrange's theorem
we obtain
\begin{eqnarray*}
  B  & \leq & \sum_{n=1}^{L} \int_{B(z-x,2^{n+1}|h|)\setminus B(z-x,2^{n}|h|)} 
                |f_h(v)|^{1+\varepsilon}\nu(dv)\\
      & \leq & \sum_{n=1}^{L} c (2^{n-1}|h|)^{(\alpha-d-1)(1+\varepsilon)}
                               \nu(B(z-x,2^{n+1}|h|)) \\
      & \leq & \sum_{n=1}^{L} c  (2^{n-1}|h|)^{(\alpha-d-1)(1+\varepsilon)}
               |z-x|^{-\alpha}\nu(B(\frac{z-x}{|z-x|},\frac{2^{n+1}|h|}{|z-x|})) \\
      & \leq & \sum_{n=1}^{L} c (2^{n-1}|h|)^{(\alpha-d-1)(1+\varepsilon)}
               (2^{n+1}|h|)^{\gamma} |z-x|^{-\alpha-\gamma}\\
      &   =  & c |h|^{(\alpha-d-1)(1+\varepsilon)+\gamma} |z-x|^{-\alpha-\gamma}
               \sum_{n=1}^{L} (2^{(\alpha-d-1)(1+\varepsilon)+\gamma})^n \\
      & \leq & c |z|^{-\alpha-\gamma}.
\end{eqnarray*}
Also by (\ref{pochodna_Green_osz}) we get
\begin{eqnarray*}
  C &  =   & \int_{B(z,1)\setminus B(z-x,1/8)} |f_h(v)|^{1+\varepsilon}\nu(dv) \\
    & \leq & \int_{B(z,1)\setminus B(z-x,1/8)}
             c\nu(dv) \leq c \nu(B(z,1)).
\end{eqnarray*}
We obtained
$$
  \int_{B(z,1)}|f_h(v)|^{1+\varepsilon} \nu(dv) = A+B+C \leq 
   c(|z|^{-\alpha-\gamma}+\nu(B(z,1))),
$$
which yields the uniform integrability of the family of functions $\{ f_h(v)\}$.
This, together with Proposition \ref{Green_regular} and (\ref{e:P}), 
yields
\begin{equation}
  D_{x_i} k_z(x) = \int_{B(z,1)}  
  D_{x_i} G_B(x,z-v)\nu(dv)\,,\quad |x|<1/2,\, |z|>1,
\end{equation}
and by (\ref{pochodna_Green_osz}), Lemma \ref{calka_gamma_miary} with $a=d-\alpha+1$ and Lemma \ref{miara_L} we get also (\ref{Poisson_pochodna_szac}) (compare the proof of Lemma \ref{Poissonszacbrzeg}).
\\
Let now $\kappa_1>0$, $\kappa_1\neq 1$ and $|\beta|=\czc{\kappa_1}=\czcd{\kappa_1}$. We note that $\kappa_1<2$ so that 
$|\beta|\in\{0,1\}$. Let $|z|> 1$, $x,y\in B(0,1/2)$ and $|y-x|<1/16$. We have
\begin{eqnarray*}
  |D^\beta k_z(y)-D^\beta k_z(x)| & \leq & \int_{B(z,1)} |D^\beta_x G_B(y,z-v)-D^\beta_x G_B(x,z-v)|\nu(dv) \\
                &   =  & \int_{B(z-x,2|y-x|)} + \int_{B(z-x,1/8)\setminus B(z-x,2|y-x|)} \\
                &      & +  \int_{B(z,1)\setminus B(z-x,1/8)} \\
                &   =  & I + II + III.
\end{eqnarray*}
By (\ref{pochodna_Green_osz}) and Lemma \ref{calka_gamma_miary} with $a=d-\alpha+|\beta|$ we have
\begin{eqnarray*}
  I & \leq & \int_{B(z-y,3|y-x|)} |D^\beta_x G_B(y,z-v)|\nu(dz) \\ 
    &      & + \int_{B(z-x,3|y-x|)}|D^\beta_x G_B(x,z-v)|\nu(dv) \\ 
    & \leq & \int_{B(z-y,3|y-x|)} c |z-y-v|^{\alpha-d-|\beta|}\nu(dv) \\ 
    &      &  + \int_{B(z-x,3|y-x|)} c |z-x-v|^{\alpha-d-|\beta|}\nu(dv) \\ 
    & \leq & c |y-x|^{\czu{\kappa_1}} |z-y|^{-\alpha-\gamma}+
             c |y-x|^{\czu{\kappa_1}} |z-x|^{-\alpha-\gamma} \\
    & \leq & c |y-x|^{\czu{\kappa_1}} |z|^{-\alpha-\gamma}.
\end{eqnarray*}
In order to estimate $II$ and $III$ we assume first $\kappa_0\not\in\N$ or $\kappa_0-\kappa_1=\alpha>1$. Let $\delta_0=\min( \kappa_0,|\beta|+1 )$. Then we have $\czcd{\delta_0}=|\beta|=\czc{\kappa_1}$, because 
$\kappa_0>\kappa_1>\czc{\kappa_1}$. 
We use (\ref{holder_Green}) in case $\delta_0=\kappa_0$, and (\ref{pochodna_Green_osz}) and
Lagrange's theorem in case $\delta_0=|\beta|+1$, to obtain for every
$r\in(0,1/2)$ and $w\in B$
\begin{equation}
  |D^\beta_x G_B(y,w)-D^\beta_x G_B(x,w)|  \leq   c r^{\alpha-d-\delta_0}|y-x|^{\delta_0-\czcd{\delta_0}},
\end{equation}
provided $x,y\in B(0,1/2)\cap B(w,r)^c$. Hence for $L=\czc{-\log_2(8|y-x|)} $ we get
\begin{eqnarray*}
  II  & \leq & \sum_{n=1}^{L} \int_{B(z-x,2^{n+1}|y-x|)\setminus B(z-x,2^{n}|y-x|)} 
               |D^\beta_x G_B(y,z-v)-D^\beta_x G_B(x,z-v)|\nu(dv)\\
      & \leq & \sum_{n=1}^{L} c (2^{n-1}|y-x|)^{\alpha-d-\delta_0} |y-x|^{\delta_0-\czcd{\delta_0}}
               \nu(B(z-x,2^{n+1}|y-x|)) \\
      & \leq & \sum_{n=1}^{L} c  (2^{n-1})^{\alpha-d-\delta_0}|y-x|^{\alpha-d-\czcd{\delta_0}}
               |z-x|^{-\alpha}\nu(B(\frac{z-x}{|z-x|},\frac{2^{n+1}|y-x|}{|z-x|})) \\
      & \leq & \sum_{n=1}^{L} c (2^{n-1})^{\alpha-d-\delta_0}
               |y-x|^{\alpha-d-\czcd{\delta_0}}(2^{n+1}|y-x|)^{\gamma} |z-x|^{-\alpha-\gamma}\\
      & \leq & c |y-x|^{\alpha+\gamma-d-\czcd{\delta_0}} |z-x|^{-\alpha-\gamma} 
               \sum_{n=1}^{\infty} (2^{\alpha-d-\delta_0+\gamma})^n \\
      & \leq & c |y-x|^{\kappa_1-\czcd{\delta_0}} |z|^{-\alpha-\gamma}\\
      &   =  & c |y-x|^{\czu{\kappa_1}} |z|^{-\alpha-\gamma},
\end{eqnarray*}
and 
\begin{eqnarray*}\label{szac_III}
  III &  =   & \int_{B(z,1)\setminus B(z-x,1/8)} 
               |D^\beta_x G_B(y,z-v)-D^\beta_x G_B(x,z-v)|\nu(dv) \\ 
      & \leq & \int_{B(z,1)\setminus B(z-x,1/8)} c (1/16)^{\alpha-d-\delta_0}
               |y-x|^{\delta_0-\czcd{\delta_0}} \nu(dv)\\ 																							
      & \leq & c |y-x|^{\delta_0-\czcd{\delta_0}} \nu(B(z,1))\\
      & \leq & c |y-x|^{\delta_0-\czcd{\delta_0}} |z|^{-\gamma} (|z|-1)^{-\alpha}.
\end{eqnarray*}
This, together with Lemma \ref{Poissonszacbrzeg} and (\ref{Poisson_pochodna_szac}), yield (\ref{Poisson_holder}) 
for $\kappa_0 \not\in\N$ or $\kappa_0-\kappa_1>1$.
If $\kappa_0\in\N$ and $\kappa_0-\kappa_1< 1$
then by (\ref{holder_Green}) we obtain
\begin{eqnarray*}
  II  & \leq & \sum_{n=1}^{L} \int_{B(z-x,2^{n+1}|y-x|)\setminus B(z-x,2^{n}|y-x|)} 
               |D^\beta_x G_B(y,z-v)-D^\beta_x G_B(x,z-v)|\nu(dv)\\
      & \leq & \sum_{n=1}^{L} c(2^{n-1}|y-x|)^{\alpha-d-|\beta|} 2^{-n+1}\log(2^{n+1})
               \nu(B(z-x,2^{n+1}|y-x|)) \\
      & \leq & \sum_{n=1}^{L} c n  (2^{n-1})^{\alpha-d-|\beta|-1}|y-x|^{\alpha-d-|\beta|}
               |z-x|^{-\alpha}\nu(B(\frac{z-x}{|z-x|},\frac{2^{n+1}|y-x|}{|z-x|})) \\
      & \leq & \sum_{n=1}^{L} c n(2^{n-1})^{\alpha-d-|\beta|-1}
               |y-x|^{\alpha-d-|\beta|}(2^{n+1}|y-x|)^{\gamma} |z-x|^{-\alpha-\gamma}\\
      & \leq & c |y-x|^{\kappa_1-|\beta|} |z-x|^{-\alpha-\gamma} \sum_{n=1}^{\infty} n(2^{-\alpha})^n \\
      & \leq & c |y-x|^{\czu{\kappa_1}} |z|^{-\alpha-\gamma},
\end{eqnarray*}
and
\begin{eqnarray*}
  III &  =   & \int_{B(z,1)\setminus B(z-x,1/8)} 
               |D^\beta_x G_B(y,z-v)-D^\beta_x G_B(x,z-v)|\nu(dv) \\
      & \leq & \int_{B(z,1)\setminus B(z-x,1/8)}
               c(1/16)^{\alpha-d-\kappa_0}\log(2+\frac{1}{16|y-x|})|y-x| \nu(dv)\\
      & \leq & c \log(2+\frac{1}{16|y-x|})|y-x| \nu(B(z,1))\\
      & \leq & c |y-x|^{\czu{\kappa_1}} |z|^{-\gamma} (|z|-1)^{-\alpha},
\end{eqnarray*}
and we get (\ref{Poisson_holder}) also in this case. \\
If $\kappa_1=1$ and $\beta=0$ then by Lemma \ref{calka_gamma_miary} with $a=d-\alpha$, (\ref{pochodna_Green_osz}) and
Lagrange's theorem we get
\begin{eqnarray*}
  I & \leq & \int_{B(z-y,3|y-x|)} |G_B(y,z-v)|\nu(dz) \\
    &      & + \int_{B(z-x,3|y-x|)}|G_B(x,z-v)|\nu(dv) \\
    & \leq & c |y-x| |z|^{-\alpha-\gamma}, \nonumber
\end{eqnarray*}
\begin{eqnarray*}
II    & \leq & c |y-x|^{\alpha+\gamma-d} |z-x|^{-\alpha-\gamma} \sum_{n=1}^{L} (2^{\alpha-d-1+\gamma})^n \\
      &  =   & c |y-x| |z-x|^{-\alpha-\gamma} L \\
      & \leq & c |y-x| |z-x|^{-\alpha-\gamma} \log(|y-x|^{-1}),
\end{eqnarray*}
and
\begin{eqnarray*}
  III &  =   & \int_{B(z,1)\setminus B(z-x,1/8)} 
               |G_B(y,z-v)-G_B(x,z-v)|\nu(dv) \\
      & \leq & c |y-x| |z|^{-\gamma} (|z|-1)^{-\alpha},
\end{eqnarray*}

which yields (\ref{Poisson_holder}) also in this case. By (\ref{Poisson_szac}), (\ref{Poisson_pochodna_szac}) and (\ref{Poisson_holder}) we get $k_z\in C^{\kappa_1}(B(0,1/2))$ for $\kappa_1\neq 1$ and 
$k_z\in C^p(B(0,1/2))$ for $\kappa_1=1$ and $p\in(0,1)$.
\qed
\section{Harmonic functions}\label{s:hf}

The H\"older continuity of harmonic functions of stable processes under assumption that
$\nu$ is absolutely continuous was first proved 
in \cite{BL}. Here we do not generally assume that $\nu$ is absolutely continuous (or $\gamma=d$).
However, if so is the case, then we obtain $\kappa_1=\alpha$ and we get the $\alpha$-H\"older continuity for $\alpha<1$ and the existence of the first
derivative of bounded harmonic functions for $\alpha>1$.

The following lemma contains continuity estimates for bounded harmonic functions. These inequalities
yield H\"older continuity.

\begin{lemat}\label{harmon_holder}
 If $\kappa_1=\gamma-(d-\alpha)>0$ then there exists a constant $c$ such that
for every function $u$ bounded on $\R$ and harmonic in $B$ we have
$$
  |u(y)-u(x)| \leq c \|u \|_\infty f(|y-x|)\,, \quad x,y\in B(0,1/2),
$$
where 
$$
f(\rho)=
\left\{
\begin{array}{lcl}
  \rho^{\kappa_1}   & \mbox{for }  & \alpha<1,\\
  \rho^{\kappa_1}\log(2/\rho)    & \mbox{for }  & \alpha=1 \mbox{  and  } \gamma<d,\\
  \rho\log^2(2/\rho)    & \mbox{for }  & \alpha=1 \mbox{  and  } \gamma=d,\\
  \rho^{\frac{2-\alpha}{\alpha}\min(\kappa_1,1)}
                & \mbox{for } & \alpha>1 \mbox{  and  } \kappa_1\neq 1, \\
  \left[\rho\log(2/\rho)\right]^{\frac{2-\alpha}{\alpha}} & \mbox{for } & \alpha>1 \mbox{  and  } \kappa_1 =1.
\end{array}\right.
$$
\end{lemat}
\dowod Let $\delta_0=\min(\kappa_1,1)$. If $\kappa_1>0$ and $\kappa_1\neq 1$ then by Lemma \ref{Poisson_regularnosc}, \ref{Poissonszacbrzeg} and Lagrange's theorem for $\alpha>1$ we obtain
\begin{eqnarray*}
  |u(y)-u(x)| & \leq & \int_{B^c} |u(z)| |P_B(y,z)-P_B(x,z)| dz \\
              & \leq & \|u \|_\infty \int_{B(0,2)^c}
                                     c|y-x|^{\delta_0}|z|^{-\gamma-\alpha} dz \\
              &      & + \|u \|_\infty \int_{B(0,2)\setminus B(0,1+|y-x|^{2\delta_0/\alpha})}
                                     c|y-x|^{\delta_0}(|z|-1)^{-\alpha} dz \\
              &      & + \|u \|_\infty \int_{B(0,1+|y-x|^{2\delta_0/\alpha})\setminus B(0,1)}
                                     c(|z|-1)^{-\alpha/2} dz \\
              &   =  & c \|u \|_\infty (|y-x|^{\delta_0}+|y-x|^{\delta_0\frac{2-\alpha}{\alpha}}) \\
              & \leq & c \|u \|_\infty |y-x|^{\delta_0\frac{2-\alpha}{\alpha}}.
\end{eqnarray*}
For $\alpha<1$ by Lemma \ref{Poisson_regularnosc} we have
\begin{eqnarray*}
  |u(y)-u(x)| & \leq & \int_{B^c} |u(z)| |P_B(y,z)-P_B(x,z)| dz \\
              & \leq & \|u \|_\infty \int_{B^c}
                                     c|y-x|^{\delta_0}|z|^{-\gamma}(|z|-1)^{-\alpha} dz \\
              & =    & c \|u \|_\infty |y-x|^{\delta_0}.
\end{eqnarray*}
Similarly for $\alpha=1$ we get
$$
  |u(y)-u(x)| \leq  c \|u \|_\infty |y-x|^{\delta_0}\log(\frac{2}{|y-x|}).
$$
Let $r=[\log(2|y-x|^{-1})|y-x|]^{2/\alpha}$. If $\kappa_1=1$ then for $\alpha> 1$ 
by Lemma \ref{Poisson_regularnosc} we obtain
\begin{eqnarray*}
  |u(y)-u(x)| & \leq & \int_{B^c} |u(z)| |P_B(y,z)-P_B(x,z)| dz \\
              & \leq & \|u \|_\infty \int_{B(0,2)^c}
                                     c\log(2|y-x|^{-1})|y-x||z|^{-\gamma-\alpha} dz \\
              &      & + \|u \|_\infty 
                       \int_{B(0,2)\setminus B(0,1+r))}
                         c\log(2|y-x|^{-1})|y-x|(|z|-1)^{-\alpha} dz \\
              &      & + \|u \|_\infty 
                         \int_{B(0,1+r)\setminus B(0,1)}
                          c(|z|-1)^{-\alpha/2} dz \\
              & \leq & c \|u \|_\infty \left[\log(2|y-x|^{-1})|y-x|+
                       [\log(2|y-x|^{-1})|y-x|]^{\frac{2-\alpha}{\alpha}}\right] \\
              & \leq & c \|u \|_\infty \left[\log(2|y-x|^{-1})|y-x|\right]^{\frac{2-\alpha}{\alpha}},
\end{eqnarray*}
and for $\alpha=1$, $\gamma=d$ we have
\begin{eqnarray*}
  |u(y)-u(x)| & \leq & \int_{B^c} |u(z)| |P_B(y,z)-P_B(x,z)| dz \\
              & \leq & \|u \|_\infty \int_{B(0,2)^c}
                                     c\log(2|y-x|^{-1})|y-x||z|^{-d-1} dz \\
              &      & + \|u \|_\infty 
                       \int_{B(0,2)\setminus B(0,1+r))}
                         c\log(2|y-x|^{-1})|y-x|(|z|-1)^{-1} dz \\
              &      & + \|u \|_\infty 
                         \int_{B(0,1+r)\setminus B(0,1)}
                          c(|z|-1)^{-1/2} dz \\
              & \leq & c \|u \|_\infty \log(2|y-x|^{-1})|y-x|
                       \log\left[(\log(2|y-x|^{-1})|y-x|)^{-1}\right] , \\
              & \leq & c \|u \|_\infty \log^2(2|y-x|^{-1})|y-x|.
\end{eqnarray*}
\qed

Note that translation invariance and the scaling property of the process 
yields that if the function $y\to u(y)$ is harmonic in $B(z,r)$ then
the function $y\to u(ry+z)$ is harmonic in $B$. Therefore we have
\begin{equation}\label{skalhar}
  |u(y)-u(x)| \leq c \|u\|_\infty f(|y-x|/r)\,,\quad x,y\in B(z,r/2),
\end{equation}
for every $z\in\R$, $r>0$ and bounded function $u$ harmonic in $B(z,r)$.

\begin{lemat}\label{Poisson_zgrubsza} If $\kappa_1=\gamma-(d-\alpha)>0$ then there exists a constant $c$ such that 
for every bounded open set $D$ 
we have
\begin{equation}
  P_D(x,z) \leq c \dist(x,D^c)^{\alpha-d} \dist (z,D)^{-\alpha},\; x\in D,\, z\in D^c.
\end{equation}
\end{lemat}
\dowod Let $x\in D$ and $z\in D^c$. By (\ref{Green_podstawa}), (\ref{e:P}) and Lemma \ref{calka_gamma_miary} we have
\begin{eqnarray*}
  P_D(x,z) &  =   & \int_{z-D} G_D(x,z-v)\nu(dv) \\
           & \leq & c \int_{z-D} |z-v-x|^{\alpha-d} \nu(dv) \\
           & \leq & c \int_{B(z-x,\dist(x,D^c)/2)} |z-v-x|^{\alpha-d} \nu(dv) \\
           &      & + c \int_{z-D} (\dist(x,D^c)/2)^{\alpha-d} \nu(dv) \\
           & \leq & c [\dist(x,D^c)]^{\gamma+\alpha-d}|z-x|^{-\alpha-\gamma} \\
           &      &   + c (\dist(x,D^c))^{\alpha-d} \nu(B(0,\dist(z,D))^c) \\
           & \leq & c \dist(x,D^c)^{\alpha-d} \dist (z,D)^{-\alpha}.
\end{eqnarray*}
\qed

\begin{lemat}\label{dpothar} If $\kappa_1=\gamma-(d-\alpha)>0$  and $\alpha>1$ then
for every $i\in\{ 1,\dots ,d\}$ the function $D_i V(x)$ is harmonic in $\R\setminus \{ 0\}$.
\end{lemat}
\dowod
Let $R>r>0$, $h\in (-r/16,r/16)$ and $x\in A=B(0,R)\setminus B(0,r)$.
Since the potential kernel $V(x)$ is harmonic in $\R\setminus \{ 0\}$ we have
\begin{displaymath}
  V(x+he_i)= E^{x+he_i}V(X_{\tau_{(A+he_i)}}) = E^x V(X_{\tau_A}+he_i),
\end{displaymath}
and
\begin{eqnarray*}
  \frac{1}{h}(V(x+he_i)-V(x)) & =  & \frac{1}{h}E^x\left[V(X_{\tau_A}+he_i)-V(X_{\tau_A})\right] \\
                              & =  & \frac{1}{h} \int_{A^c} \left[V(z+he_i)-V(z)\right] \omega^x_A(dz).
\end{eqnarray*}
 Let
$$
  f_h(z)=\frac{ V(z+he_i)-V(z) }{h}.
$$
We will prove that the family of functions 
$\{ f_h(z), |h|<r/16 \}$ is uniformly integrable with respect to the measure
$\omega_A^x$.
To this end we will estimate the integral
\begin{eqnarray*}
  \int_{A^c}|f_h(z)|^{1+\varepsilon} \omega_A^x(dz) & = & \int_{B(0,2|h|)} +
   																							 \int_{B(0,r/8)\setminus B(0,2|h|)}
   																							+ \int_{A^c \setminus B(0,r/8)} \\
   																							& = & I+II+III,
\end{eqnarray*}
where $0<\varepsilon<\frac{d}{d-\alpha+1}-1$.
By Lemma \ref{Poisson_zgrubsza} we have
\begin{eqnarray*}
  I & \leq & \frac{2^{1+\varepsilon}}{|h|^{1+\varepsilon}} 
               \left[\int_{B(-he_i,3|h|)} (V(z+he_i))^{1+\varepsilon}\omega_A^x(dz)  
                + \int_{B(0,3|h|)} (V(z))^{1+\varepsilon}\omega_A^x(dz)\right] \\
    & \leq & \frac{2^{1+\varepsilon}}{|h|^{1+\varepsilon}} [ \int_{B(-he_i,3|h|)} 
                                  c|z+he_i|^{(\alpha-d)(1+\varepsilon)}P_A(x,z)dz \\
    &      &  + \int_{B(0,3|h|)} c |z|^{(\alpha-d)(1+\varepsilon)}P_A(x,z)dz] \\
    & \leq & \frac{c}{|h|^{1+\varepsilon}} |h|^{d+(\alpha-d)(1+\varepsilon)} 
              \dist(x,A^c)^{\alpha-d} (3r/4)^{-\alpha}  \\
    & \leq & c \dist(x,A^c)^{\alpha-d} (3r/4)^{-\alpha}.
\end{eqnarray*}
By (\ref{pochodna_osz}) and Lagranger's theorem
for $L=\czc{\log_2(\frac{r}{8|h|})}$ we obtain
\begin{eqnarray*}
  II  & \leq & \sum_{n=1}^{L} \int_{B(0,2^{n+1}|h|)\setminus B(0,2^{n}|h|)} 
                |f_h(z)|^{1+\varepsilon}\omega^x_A(dz)\\
      & \leq & \sum_{n=1}^{L} c (2^{n-1}|h|)^{(\alpha-d-1)(1+\varepsilon)}
                               \omega^x_A(B(0,2^{n+1}|h|)) \\
      & \leq & \sum_{n=1}^{L} c  (2^{n-1}|h|)^{(\alpha-d-1)(1+\varepsilon)}
               \dist(x,A^c)^{\alpha-d} (3r/4)^{-\alpha} (2^{n+1}|h|)^d \\
      &   =  & c |h|^{(\alpha-d-1)(1+\varepsilon)+d} \dist(x,A^c)^{\alpha-d} (3r/4)^{-\alpha}
               \sum_{n=1}^{L} (2^{(\alpha-d-1)(1+\varepsilon)+d})^n \\
      & \leq & c \dist(x,A^c)^{\alpha-d} r^{(\alpha-d-1)(1+\varepsilon)+d-\alpha},
\end{eqnarray*}
\begin{eqnarray*}
  III &  =   & \int_{A^c\setminus B(0,r/8)} |f_h(z)|^{1+\varepsilon}\omega_A^x(dz) \\
    & \leq & \int_{A^c\setminus B(0,r/8)}
             cr^{(\alpha-d-1)(1+\varepsilon)}\omega_A^x(dz)\\
      & \leq & cr^{(\alpha-d-1)(1+\varepsilon)} \omega_A^x(A^c)=cr^{(\alpha-d-1)(1+\varepsilon)}.
\end{eqnarray*}
We obtained
$$
  \int_{A^c}|f_h(z)|^{1+\varepsilon} \omega_A^x(dz) = I+II+III \leq c(r,R,x)
$$
which clearly yields the uniform integrability of the family of functions $\{ f_h(v)\}$. We obtain
\begin{eqnarray*}
  D_i V(x)  & = & \lim_{h\to 0} \frac{V(x+he_i)-V(x)}{h} \\
            & = & \lim_{h\to 0} \int_{A^c} \frac{V(z+he_i)-V(z)}{h} \omega^x_A(dz) \\
            & = & \int_{A^c} D_i V(z) \omega^x_A(dz) \\
            & = & E^x D_i V(X_{\tau_A}).
\end{eqnarray*}
This yields regular harmonicity of $D_iV$ in $A=B(0,R)\setminus B(0,r)$ for every $R>r>0$ and harmonicity in $\R\setminus \{ 0\}$.
\qed

\begin{wniosek}\label{Green_harm}
 If $\kappa_1=\gamma-(d-\alpha)>0$ and $\alpha>1$ then for every $i\in \{ 1,\dots ,d\}$ and $x\in B$ the function 
 $g(v)=D_{x_i} G_B(x,v) $ is harmonic in $B\setminus\{x \}$ and continuous in $\R\setminus \{ x\}$.
\end{wniosek}
\dowod By symmetry of the Green function, (\ref{wzornaGreena}) and Lemma \ref{potencjalmiary} we have
\begin{equation}
  D_{x_i}G_B(x,v)=
-D_i V(v-x)-E^v D_i V(x-X_{\tau_B})
\end{equation}
which is harmonic in $B\setminus\{ x\}$ by Lemma \ref{dpothar}.
The function $v\to D_i V(v-x)$ is continuous in $\R\setminus \{x \}$ by Theorem 
\ref{jadropotencjalu}. The function
$$
  h(v)=E^v D_i V(x-X_{\tau_B})\,,\quad v\in B,
$$
is bounded and harmonic in $B(v_0,1-|v_0|)$ for every $v_0\in B$ hence from (\ref{skalhar})  we obtain the continuity of $h(v)$ at $v_0$ and the continuity of $v\to D_{x_i}G_B(x,v)$ in $B$. $D_{x_i}G_B(x,v)=0$ for $v\in B^c$ and we have also $\lim\limits_{v\to w} D_{x_i}G_B(x,v)=0$ for every
$x\in B(0,1)$ and every point $w\in \sfera$ because 
the measures $\omega^v_{B}$ weakly converge to $\delta_{w}$ (see, e.g., \cite[Theorem 1.23]{ChZ}).
\qed

We will employ the operator
 $$
  \Dynkinr{\phi}{x} = \frac{E^x\phi(X_{\tau_{B(x,r)}}) - \phi(x) } {E^x \tau_{B(x,r)}}\,,
 $$
  whenever the expression is well defined for given $\phi$, $r>0$
  and $x$. 
Clearly, if $h$ is harmonic in $D$, $x\in D$, and
$r<\dist(x,D^c)$, then $\Dynkinr{h}{x}=0$. We note that
$$
  \Dynkin{\phi}{x}=\lim_{r\downarrow 0} \Dynkinr{\phi}{x}
$$
is the Dynkin characteristic operator, which was used
in \cite{BS} in a similar way.

We record the following observation (maximum principle).
\begin{lemat}\label{zasadaminimum}
  If there is $r>0$ such that $\Dynkinr{h}{x} > 0$ then $h(x)<\sup\limits_{y\in\R} h(y)$.
\end{lemat}

The proof of the following Lemma is based on ideas of \cite[Lemma 16]{BS}.

\begin{lemat}\label{Greenpoch<s}
  If $\kappa_1=\gamma-(d-\alpha)>1$ then there exists $c$ such
  that
  \begin{equation}\label{Greenpochszacbrzeg}
   |D_{x_i}G_B(x,v)| < c (1-|v|)^{\alpha/2}\,,\quad |x|<1/2\,,\;3/4<|v|<1\,.
  \end{equation}
\end{lemat}
\dowod 
By the strong Markov property we have
\begin{eqnarray*}
  s(v)
  &   =   & E^v \tau_B = E^v(\tau_A+\tau_{B(0,1)} \circ \theta_{\tau_A})
     =    E^v\tau_A + E^v E^{X_{\tau_A}} \tau_{B(0,1)} \\
  &   =   & E^v \tau_A + E^v s(X_{\tau_A})\,,\quad v\in \R\,,\quad A\subset {B(0,1)},
\end{eqnarray*}
which yields $\Dynkinr{s}{v}=-1$ for $v\in {B(0,1)}$ and $r<1-|v|$.

For $n\in \N$ and  $x\in B(0,1/2)$ we let 
$g(v) =D_{x_i}G_B(x,v)$, $g_{1,n}(v) =\min(g(v),n)$ and $g_{2,n}(v)=\max (g(v),-n) $. For $v \in B(x,1/8)^c$ we have that
$|D_{x_i}G_B(x,v)| \leq c_1|v-x|^{\alpha-d-1}$ hence
$g_{1,n}(v)=g_{2,n}(v)=g(v)$ provided $n \geq c_1 8^{d-\alpha+1}$. It follows from
Corollary \ref{Green_harm} that $g(v)$ is harmonic on $B\setminus\{x\}$. Hence by
scaling property, (\ref{IW}), Lemma \ref{calka_gamma_miary} and \ref{miara_L}
we obtain  that for $v\in {B}\setminus B(0,3/4)$ and $r<\min(1-|v|,1/16)$ it holds
\begin{eqnarray*}
  \Dynkinr{g_{1,n}}{v} 
 &  =    &    \Dynkinr{(g_{1,n}-g)}{v} \\
 &  =   & \frac{1}{E^0\tau_{B}}
           \int\limits_{B} G(0,w)
           \int(g_{1,n}-g)(v+rw+z) \nu(dz)dw \\
  & \geq & \frac{-c_2}{s(0)} \int\limits_{B} G(0,w)
           \int\limits_{B(x-v-rw,1/8)}|x-v-rw-z|^{\alpha-d-1} \nu(dz)dw \\
  & \geq & -c_3\, ,
\end{eqnarray*}
Similarly
\begin{eqnarray*}
  \Dynkinr{g_{2,n}}{v} 
 &  =    &    \Dynkinr{(g_{2,n}-g)}{v} \\
 & \leq & c_4\, ,
\end{eqnarray*}
If $a>c_3$ and $a>c_4$ then
$$
\Dynkinr{(a s - g_{1,n})}{v} = -a - \Dynkinr{g_{1,n}}{v} \leq -a + c_3 < 0,
$$
and
$$
\Dynkinr{(a s + g_{2,n})}{v} = -a + \Dynkinr{g_{2,n}}{v} \leq -a + c_4 < 0.
$$ 
By scaling
\begin{eqnarray}\label{sdodatnie}
s(v) 
& \geq & E^v \tau_{B(v,1-|v|)}=(1-|v|)^\alpha E^0\tau_{B(0,1)} \\
& \geq & 4^{-\alpha} E^0\tau_{B}\,, \quad |v|<3/4\,. \nonumber
\end{eqnarray}
Since $g_{1,n}(v) \leq n$ and $g_{2,n}(v) \geq -n$, we see that $as(v) - g_{1,n}(v) > 0$ and $as(v)+g_{2,n}(v)> 0$ for $v\in B(0,3/4)$ provided $a > n/(4^{-\alpha}E^0\tau_{B} )$.

Let $a_0= \max[c_3,c_4,n/(4^{-\alpha}E^0\tau_{B(0,1)})] + 1$ and $h_1(v) =
a_0 s(v) - g_{1,n}(v)$, $h_2(v) =
a_0 s(v) + g_{2,n}(v)$.  We have $h_1(v)\geq 0 $ and $h_2(v)\geq 0$ for $v\in
\overline{B(0,3/4)}$, $h_1(v)=h_2(v)=0$ for $v\in {B(0,1)}^c$ and $\Dynkinr{h_1}{v} < 0$, $\Dynkinr{h_2}{v} < 0$
for $v\in {B(0,1)}\setminus B(0,3/4)$, $r<\min(1-|v|,1/16)$.
Lemma \ref{zasadaminimum} and 
continuity of $h_1$ and $h_2$ in $\R\setminus \{ x\}$ yields $h_1(v)\geq 0 $ and $h_2(v)\geq 0$ in ${B(0,1)}$.  Since $g_{1,n}=g_{2,n}=g$ on
$B(0,3/4)^c$, the Lemma follows from Lemma \ref{ET_est}.
\qed
\vspace{2mm}

Using (\ref{Greenpochszacbrzeg}) we obtain the following estimate. We omit the proof since it is analogous
to the proof of Lemma \ref{Poissonszacbrzeg}
\begin{lemat}\label{l:DP}
  If $\kappa_1=\gamma-(d-\alpha)>1$ then there exists a constant $c$ such that for every $i\in\{1,\dots,d\}$ we have
  \begin{equation}\label{Poissonpochszac}
    |D_{x_i} P_B(x,z)| \leq c |z|^{-\gamma}(|z|^2-1)^{-\alpha/2}\,,\quad |x|<1/2\,,\,|z|>1.
  \end{equation}
\end{lemat}

Now we can prove our main result.\\
\dowod[ of Theorem \ref{t:harmoniczne}] We assume first $\kappa_1=\gamma-(d-\alpha)>1$. By (\ref{Poissonpochszac}) we obtain
$$
  D_i u(x) = \int_{B^c} u(z) D_{x_i}P_B(x,z) dz\, , \quad x\in B(0,1/2),
$$
and by Lemma \ref{Poisson_regularnosc} and (\ref{Poissonpochszac}) we get
\begin{eqnarray}\label{harm_pholder} \nonumber
  |D_i u(y)-D_i u(x)| & \leq & \int_{B^c} |u(z)| |D_{x_i} P_B(y,z)-D_{x_i} P_B(x,z)| dz \\ \nonumber
              & \leq & \|u \|_\infty \int_{B(0,2)^c}
                                     c|y-x|^{\czu{\kappa_1}} |z|^{-\gamma-\alpha} dz \\ \nonumber
              &      & + \|u \|_\infty \int_{B(0,2)\setminus B(0,1+|y-x|^{2\czu{\kappa_1}/\alpha})}
                                     c|y-x|^{\czu{\kappa_1}}(|z|-1)^{-\alpha} dz \\
              &      & + \|u \|_\infty \int_{B(0,1+|y-x|^{2\czu{\kappa_1}/\alpha})\setminus B(0,1)}
                                     c(|z|-1)^{-\alpha/2} dz \\ \nonumber
              &   =  & c \|u \|_\infty (|y-x|^{\czu{\kappa_1}}+|y-x|^{\czu{\kappa_1}\frac{2-\alpha}{\alpha}}) \\ \nonumber
              & \leq & c \|u \|_\infty |y-x|^{\czu{\kappa_1}\frac{2-\alpha}{\alpha}}\, , \quad x,y\in B(0,1/2).
\end{eqnarray}
This and Lemma \ref{harmon_holder} yield for every $\kappa_1>0$ some constant $\delta\in(0,1)$ such that $u\in C^{\czcd{\kappa_1}+\delta}(B(0,1/2))$.
 Let now $x_0\in B$. The function $u(x)$ is harmonic in $B(x_0,1-|x_0|)$ therefore the function $h(x)=u((1-|x_0|)x+x_0)$ is
harmonic in $B$ and for $|\beta|=\czcd{\kappa_1}$ we get
\begin{eqnarray*}
  |D^\beta u(y)-D^\beta u(x_0)|  &    =  & (1-x_0)^{-|\beta|}|(D^\beta h)(\frac{y-x_0}{1-|x_0|})-(D^\beta h)(0)| \\
                                 & \leq  & c (1-|x_0|)^{-|\beta|-\delta }
                                           |y-x_0|^{\delta},\quad y\in B(x_0,\frac{1-|x_0|}{2}),
\end{eqnarray*}
which clearly yields $u\in C^{\czcd{\kappa_1}+\delta}_{\loc}(B)$ and (\ref{e:hh}).
\qed\\

{\bf Acknowledgments.} The author is grateful to Prof. R. L. Schilling  for useful discussions and his hospitality during the stay in
Marburg and to Prof. K. Bogdan
for discussions and suggestions on the paper. \\


\begin{thebibliography}{2}
\bibitem{BK} R. F. Bass, M. Kassmann, {\it H\"older continuity of harmonic functions with
respect to operators of variable orders,} Comm. Partial Differential Equations 30 (2005), 1249-1259.
\bibitem{BL} R. F. Bass, D. A. Levin, {\it Harnack inequalities
    for jump processes,} Potential Anal. 17(4)(2002), 375-388.
\bibitem{Br} J. Bertoin, {\it L\' evy processes,} Cambridge
  University Press, Cambridge, 1996.
\bibitem{BGB} R. M. Blumenthal and R. K. Getoor, {\it Markov Processes
    and Potential Theory,} Pure Appl. Math., Academic Press Inc., New
  York 1968.
\bibitem{BBsm1999}
K. Bogdan, T. Byczkowski, {\it Potential theory for the
  {$\alpha$}-stable {S}chr\"odinger operator on bounded {L}ipschitz domains,}
Studia Math. 133 (1999) no.~1, 53--92.

\bibitem{BD} K. Bogdan, B. Dyda, {\it Relative Fatou theorem for harmonic functions of rotation invariant stable processes in smooth domains,} Stud. Math. 157, No. 1, 83-96 (2003).

\bibitem{BJ} K. Bogdan, T. Jakubowski, {\it Estimates of heat kernel of 
fractional Laplacian perturbed by gradient operators}, Comm. Math. Phys. 271 (1) 2007, 179--198.

\bibitem{BKN} K. Bogdan, T. Kulczycki, A. Nowak, {\it Gradient estimates for harmonic and $q$-harmonic functions of symmetric stable processes,} Ill. J. Math. 46, No. 2, 541-556 (2002).

\bibitem{BKK} K. Bogdan, T. Kulczycki, M. Kwa\'snicki {\it Estimates and structure of $\alpha$ - harmonic functions,} Probab. Theory Relat. Fields 2007.

\bibitem{BSS1} K. Bogdan, A. St\'os, P. Sztonyk, {\it Potential theory
    for L\'evy stable processes,} Bull. Polish. Acad. Sci. Math. 50(3)
  2002, 361--372.

\bibitem{BSS2} K. Bogdan, A. St\'os, P. Sztonyk, {\it
    Harnack inequality for symmetric stable processes on fractals,} C.
  R. Math. Acad.  Sci. Paris 335 (1) (2002), 59--63.

\bibitem{BS0} K. Bogdan, P. Sztonyk, {\it Harnack's  inequality for
    stable L{\'e}vy processes,} Potential Analysis 22 (2) (2005), 133--150.

\bibitem{BS} K. Bogdan, P. Sztonyk, {\it Estimates of potential kernel and
Harnack's  inequality for anisotropic fractional Laplacian,} to appear in Studia Math. (\url{http://arxiv.org/PS_cache/math/pdf/0507/0507579.pdf}).

\bibitem{Ch} K. L. Chung, 
{\it Lectures from Markov processes to Brownian motion}, 
Springer-Verlag, New York-Berlin, 1982. 

\bibitem{Ch1} K. L. Chung, {\it Doubly-Feller process with multiplicative functional},
Seminar on stochastic processes, 1985 (Gainesville, Fla., 1985), 63--78, Progr. Probab. Statist., 12, 
Birkh\"auser Boston, Boston, MA, 1986. 

\bibitem{ChZ} K. L. Chung, Z. Zhao, {\it From Brownian motion to
   Schr\"odinger's equation,} Springer - Verlag, New York, 1995.
\bibitem{EtK} S.N. Ethier,Th. G. Kurtz, {\it Markov processes - charakterization and convergence,} 
Wiley Series in Probability and Mathematical Statistics, John Wiley, 
New York - Chicester - Brisbane - Toronto - Singapore (1986).

\bibitem{DG} E. De Giorgi, {\it Sulla differenziabilit\'a e l'analiticit\'a delle estremali degli 
integrali multipli regolari,} Mem. Accad. Sci. Torino Cl. Sci. Fis. Mat. Natur. (3), 
3 (1957), 25-43.

\bibitem{H} M. Hardy, {\it Combinatorics of Partial Derivatives}, Electronic Journal of Combinatorics, 13 (2006).
\bibitem{Hi} S. Hiraba {\it Asymptotic estimates for densities of
 multi-dimensional stable distributions},  Tsukuba J. Math.  27
 (2003),  no. 2, 261--287.
 
\bibitem{HY} B. Hu, H. -M. Yin, {\it The DeGiorgi-Nash-Moser type of estimate for parabolic
Volterra integrodifferentail equations}, Pacific J. Math., 178(2) (1997), 265-277.
 
\bibitem{HK} R. Husseini, M. Kassmann, {\it Jump processes, $\cal{L}$-harmonic functions, continuity estimates and
the Feller property,} 2006, preprint.
\bibitem{IW} N. Ikeda, S. Watanabe, {\it On some relations between the
    harmonic measure and the L\'evy measure for a certain class of
    Markov processes,} J. Math. Kyoto Univ. 2-1 (1962), 79-95.
\bibitem{Jc1}
N. Jacob, {\it Pseudo differential operators and Markov
  processes. Vol. I. Fourier analysis and semigroups},
Imp. Coll. Press, London, 2001.

\bibitem{Jc2}
 N. Jacob,
 {\it Pseudo-Differential Operators and Markov Processes,
 Vol. II : Generators and Their Potential Theory,}
 Imperial College Press, London, 2002.
 
\bibitem{Jk} T. Jakubowski, {\it The estimates for the Green function
    in Lipschitz domains for the symmetric stable processes,} Probab. Math. Statist.  22  
(2002), 419--441.

\bibitem{Kas} M. Kassmann, {\it On regularity for Beurling-Deny type Dirichlet forms}, Potential Analysis 19 (2003), 69-87. 

\bibitem{K1}
T{.} Kulczycki,
\emph{Properties of Green function of symmetric stable process},
Probab{.} Math{.} Statist{.} 17(2) (1997), 339--364.

\bibitem{Lewand} M. Lewandowski {\it Point regularity of $p$-stable density in ${\cal R}^d$ and Fisher information,} Probab. Math. Stat. 19, No.2, 375-388 (1999).

\bibitem{Mich2} K. Michalik, {\it Sharp estimates of the Green function, the Poisson kernel and the Martin kernel of cones for symmetric stable processes,} Hiroshima Math. J. 36, No. 1 (2006), 1-21.

\bibitem{Mich1} K. Michalik, M. Ryznar, {\it Relative Fatou theorem for $\alpha$-harmonic functions in Lipschitz domains,} Illinois J. Math. 48 (2004), 977-998.

\bibitem{MP} R. Mikulevi\v{c}ius, H. Pragarauskas, {\it On H\"older continuity
of solutions of certain integro-differential equations}, 
Ann. Acad. Sci. Fenn. Ser. A I Math., 13(2) (1988), 231-238.

\bibitem{Millar} P. W. Millar, {\it First passage distributions of
    processes with independent increments,} Ann. Probab., 3 (1975) no.
  2, 215-233.
  
\bibitem{Mos} J. Moser, {\it On Harnack's theorem for elliptic differential equations,} 
Comm. Pure Appl. Math., 14 (1961), 577-591.

\bibitem{Nash} J. Nash, {\it Continuity of solutions of parabolic and elliptic equations,} 
Amer. J. Math., 80 (1958), 931-954. 

\bibitem{Pi} J. Picard, 
{\it Density in small time at accessible points for jump processes,} 
Stochastic Process. Appl. 67 (1997), no.~2, 251--279.

\bibitem{Rudin} W. Rudin, {\it Functional Analysis}, 2nd ed., International Series in Pure and Applied Mathematics. New York, NY: McGraw-Hill. xviii, 424 p. (1991). 
Stochastic Process. Appl. 67 (1997), no.~2, 251--279.

\bibitem{Sato} K.-I. Sato, {\it L\' evy Processes and Infinitely Divisible 
             Distributions}, Cambridge University Press, 1999.
             
\bibitem{SU} R. Schilling, T. Uemura, {\it Dirichlet forms generated by pseudo 
differential operators: on the Feller property of the associated stochastic process,} to appear in Tohoku Math. J.

\bibitem{SV} 
R. Song, Z. Vondra\v{c}ek, {\it Harnack inequality
    for some classes of Markov processes}, Math. Z.  246  (2004),
  no. 1-2, 177--202.
\bibitem{Sztonyk} P. Sztonyk, {\it On harmonic measures for L\'evy
    processes,} Prob. Math. Statist. 20(2) (2000), 383-390.
\bibitem{Taylor} S. J. Taylor, {\it Sample path properties of a
    transient stable process,} J. Math. Mech. 16 (1967), 1229-1246.
\bibitem{Trib} H. Tribel, {\it Theory of function spaces,} Monographs in Mathematics, Vol. 78. Basel-Boston-Stuttgart: Birkh\"auser Verlag, (1983). 
\bibitem{W} T. Watanabe {\it Asymptotic estimates of
multi-dimensional stable densities and their applications}, Trans. Amer. Math. Soc. 359 (2007), 2851-2879.
\bibitem{Wu}
J.-M. Wu {\it Harmonic measures for symmetric stable processes},
Studia Math.  149  (2002),  no. 3, 281--293.
\end{thebibliography}
\end{document}